\documentclass[11pt, letterpaper]{amsart}
\usepackage[colorlinks]{hyperref}
\AtBeginDocument{
  \hypersetup{
    linkcolor=blue,
    citecolor=green,
  }
}
\usepackage{cleveref}
\crefname{equation}{}{}

\newtheorem{theorem}{Theorem}[section]
\newtheorem{lemma}[theorem]{Lemma}
\newtheorem{corollary}[theorem]{Corollary}

\theoremstyle{definition}
\newtheorem{definition}[theorem]{Definition}

\newtheorem{notation}[theorem]{Notation}

\theoremstyle{remark}
\newtheorem{remark}[theorem]{Remark}

\numberwithin{equation}{section}

\begin{document}

\title[Boundary pointwise regularity and applications]{Boundary pointwise regularity and applications to the regularity of free boundaries}

\author{Yuanyuan Lian}
\address{School of Mathematical Sciences, Shanghai Jiao Tong University, Shanghai, China}
\email{lianyuanyuan@sjtu.edu.cn; lianyuanyuan.hthk@gmail.com}

\author{Kai Zhang}
\address{School of Mathematical Sciences, Shanghai Jiao Tong University, Shanghai, China}
\email{zhangkaizfz@gmail.com}
\thanks{This research is supported by the China Postdoctoral Science Foundation (Grant No.
2021M692086), the National Natural Science Foundation of China (Grant No. 12031012 and 11831003) and  the Institute of Modern Analysis-A Frontier Research Center of Shanghai.}


\date{}



\begin{abstract}
In this paper, we develop a series of boundary pointwise regularity for Dirichlet problems and oblique derivative problems. As applications, we give direct and simple proofs of the higher regularity of the free boundaries in obstacle-type problems and one phase problems.
\end{abstract}

\maketitle

\section{Introduction}\label{S1}
In this paper, we prove some new pointwise boundary regularity for Dirichlet problems:
\begin{equation}\label{e.Dirichlet}
 \left\{ \begin{aligned}
&\Delta u=f ~~&&\mbox{in}~~\Omega\cap B_1;\\
&u=g ~~&&\mbox{on}~~\partial \Omega\cap B_1
 \end{aligned}\right.
\end{equation}
and oblique derivative problems:
\begin{equation}\label{e.oblique}
 \left\{ \begin{aligned}
&\Delta u=f ~~&&\mbox{in}~~\Omega\cap B_1;\\
&\beta\cdot Du=g ~~&&\mbox{on}~~\partial \Omega\cap B_1,
 \end{aligned}\right.
\end{equation}
where $\Omega\subset \mathbb{R}^n$ is a bounded domain and $0\in \partial \Omega$. The pointwise regularity shows a clear and deep relation between the regularity of solutions and the regularity of prescribed data. It can be tracked at least to the work of Caffarelli \cite{MR1005611} for the interior pointwise regularity of fully nonlinear elliptic equations. Various pointwise regularity have been developed by many researchers since then, such as boundary regularity (\cite{MR4088470}, \cite{MR3246039}), regularity for equations with lower terms (\cite{MR2334822}, \cite{MR3246039}, \cite{lian2020pointwise}, \cite{MR3980853}, \cite{MR2592289}), regularity for parabolic equations (\cite{MR1135923,MR1139064,MR1151267}) and regularity for the Monge-Amp\`{e}re equation (\cite{MR2983006}) etc.

In this paper, we develop a series of boundary pointwise regularity for Dirichlet problems and oblique derivative problems. We show that if the derivatives of $u$ vanish, $u$ possesses higher regularity than the usual. This was first observed in \cite{MR4088470} and we give a complete treatment for the Dirichlet problems \cref{e.Dirichlet} and the oblique derivative problems \cref{e.oblique} here.

As applications of these pointwise regularity, we prove the higher regularity of free boundaries in obstacle-type problems and one phase problems without using the partial hodograph-Legendre transformation (see \cite{MR0478079}), which is a standard method up to now. We clarify the idea briefly. Take the Dirichlet problem \cref{e.Dirichlet} for instance. It is well-known that if $\partial \Omega\in C^{k,\alpha}$ ($k\geq 1$), $u\in C^{k,\alpha}$. On the other hand, the regularity of $u$ may lead to the regularity of $\partial \Omega$ since $\varphi_i=u_i/u_n$ ($1\leq i< n$), where $\varphi$ is the representation function of $\partial \Omega$. If a problem is an overdetermined problem, i.e., we have more conditions on $u$ on the boundary, we may have higher regularity for $u$ and then higher regularity for $\partial \Omega$ and so forth. Eventually, $u$ and $\partial \Omega$ are infinite smooth.

Before stating our main results, we introduce some notations for pointwise regularity. The first is the pointwise characterization of a function in H\"{o}lder spaces, which is well-known now.
\begin{definition}\label{d-f}
Let $\Omega\subset \mathbb{R}^n$ be a bounded set (may be not a domain) and $f:\Omega\rightarrow \mathbb{R}$ be a function. We say that $f$ is $C^{k,\alpha}$ ($k\geq 0, 0<\alpha\leq 1$) at $x_0\in \Omega$ or $f\in C^{k, \alpha}(x_0)$ if there exist constants $K,r_0>0$ and a polynomial $P\in \mathcal{P}_k$ (i.e., degree less than or equal to $k$) such that
\begin{equation}\label{m-holder}
  |f(x)-P(x)|\leq K|x-x_0|^{k+\alpha},~~\forall~x\in \Omega\cap B_{r_0}(x_0).
\end{equation}
Then define $D^if(x_0)=D^iP(x_0)$ ($1\leq i\leq k$),
\begin{equation*}
\begin{aligned}
&[f]_{C^{k,\alpha}(x_0)}=\min \left\{K\big | \cref{m-holder} ~\mbox{holds with}~P~\mbox{and}~K\right\}, \\
&\|f\|_{C^{k,\alpha}(x_0)}=\|P\|+[f]_{C^{k,\alpha}(x_0)}.
\end{aligned}
\end{equation*}

If $f\in C^{k, \alpha}(x)$ for any $x\in \Omega$ with the same $r_0$ and
\begin{equation*}
  \|f\|_{C^{k,\alpha}(\bar{\Omega})}= \sup_{x\in \Omega} \|f\|_{C^{k}(x)}+\sup_{x\in \Omega} [f]_{C^{k, \alpha}(x)}<+\infty,
\end{equation*}
we say that $f\in C^{k,\alpha}(\bar{\Omega})$.

In addition, we say that $f$ is $C^{-1,\alpha}$ at $x_0$ or $f\in C^{-1,\alpha}(x_0)$ if there exist constants $K,r_0>0$ such that
\begin{equation}\label{e.c-1}
\|f\|_{L^n(\bar{\Omega}\cap B_r(x_0) )}\leq Kr^{\alpha-1}, ~\forall ~0<r<r_0.
\end{equation}
Then define
\begin{equation*}
\|f\|_{C^{-1,\alpha}(x_0)}=\min \left\{K \big | \cref{e.c-1} ~\mbox{holds with}~K\right\}.
\end{equation*}

If $f\in C^{-1, \alpha}(x)$ for any $x\in \Omega$ with the same $r_0$ and
\begin{equation*}
  \|f\|_{C^{-1,\alpha}(\bar{\Omega})}:= \sup_{x\in \Omega} \|f\|_{C^{-1,\alpha}(x)}<+\infty,
\end{equation*}
we say that $f\in C^{-1,\alpha}(\bar{\Omega})$.
\end{definition}

\begin{remark}\label{re1.1}
If $\Omega$ is a smooth domain (e.g. a Lipschitz domain), the definition of $f\in C^{k,\alpha}(\bar{\Omega})$ ($k\geq 0$) is equivalent to the classical definition.
\end{remark}

The next is a pointwise characterization of the smoothness of a domain's boundary. This definition is similar to \Cref{d-f}. That is, both definitions use polynomials to describe the smoothness. It was first introduced in \cite{MR4088470}.

\begin{definition}\label{d-domain}
Let $\Omega$ be a bounded domain, $\Gamma\subset \partial \Omega$ be relatively open and $x_0\in \Gamma$. We say that $\Gamma$ is $C^{k,\alpha}$ ($k\geq 0, 0<\alpha\leq 1$) at $x_0$ or $\Gamma\in C^{k,\alpha}(x_0)$ if there exist constants $K,r_0>0$, a coordinate system $\{x_1,...,x_n \}$ (isometric to the original coordinate system) and a polynomial $P\in \mathcal{P}_k$ with $P(0)=0$ and $DP(0)=0$ (if $k\geq 1$) such that $x_0=0$ in this coordinate system,
\begin{equation}\label{e-re}
B_{r_0} \cap \{(x',x_n)\big |x_n>P(x')+K|x'|^{k+\alpha}\} \subset B_{r_0}\cap \Omega
\end{equation}
and
\begin{equation}\label{e-re2}
B_{r_0} \cap \{(x',x_n)\big |x_n<P(x')-K|x'|^{k+\alpha}\} \subset B_{r_0}\cap \Omega^c.
\end{equation}
Then, define
\begin{equation*}
[\Gamma]_{C^{k,\alpha}(x_0)}=\min \left\{K\big | \cref{e-re} ~\mbox{and}~ \cref{e-re2}~\mbox{hold with}~P~\mbox{and}~K\right\}
\end{equation*}
and
\begin{equation*}
\|\Gamma\|_{C^{k,\alpha}(x_0)}=\|P\|+[\partial \Omega]_{C^{k,\alpha}(x_0)}.
\end{equation*}

If $\Gamma\in C^{k, \alpha}(x)$ for any $x\in \Gamma$ with the same $r_0$ and
\begin{equation*}
  \|\Gamma\|_{C^{k,\alpha}}:= \sup_{x\in \Gamma} \|\Gamma\|_{C^{k}(x)}+\sup_{x\in \Gamma}~[\Gamma]_{C^{k, \alpha}(x)}<+\infty,
\end{equation*}
we say that $\bar\Gamma\in C^{k,\alpha}$. If $\bar\Gamma'\in C^{k,\alpha}$ for any $\Gamma'\subset\subset\Gamma$, we say that $\Gamma\in C^{k,\alpha}$. If $\Gamma\in C^{k,\alpha}$ for any $k\geq1$ and $0<\alpha\leq 1$, we say that $\Gamma\in C^{\infty}$.
\end{definition}

\begin{remark}\label{re1.13}
The \cref{e-re} and \cref{e-re2} means that the difference between $\partial \Omega\cap B_{r_0}$ and $P$ is controlled by $K|x'|^{k+\alpha}$. Hence, this is an analogue of the pointwise $C^{k,\alpha}$ for a function. One benefit of \Cref{d-domain} is that $\partial \Omega\cap B_{r_0}$ doesn't need to be the graph of some function. It could be rather complicated.
\end{remark}

\begin{remark}\label{re1.2}
Throughout this paper, if we say that $f\in C^{k,\alpha}(x_0)$ ($\Gamma\in C^{k,\alpha}(x_0)$), we use $P_f$ ($P_{\Omega}$) to denote the corresponding polynomial in \Cref{d-f} (\Cref{d-domain}).
\end{remark}

\begin{remark}
We always assume that $0\in\partial \Omega$ and study the pointwise regularity at $0$ for \cref{e.Dirichlet} and \cref{e.oblique}. In addition, if we use \Cref{d-f} and \Cref{d-domain} at $0$, we always assume that $r_0=1$, and \cref{e-re} and \cref{e-re2} hold if we say $\partial \Omega \cap B_1\in C^{k,\alpha}(0)$.

We also use the following notation to describe the oscillation of $\partial \Omega$ near $0$. For $r>0$, define
\begin{equation*}
\underset{B_{r}}{\mathrm{osc}}~\partial\Omega= \underset{x\in \partial \Omega\cap B_r}{\sup} x_n -\underset{x\in \partial \Omega\cap B_r}{\inf} x_n.
\end{equation*}
\end{remark}
~\\

Now, we state our main results.
\begin{theorem}\label{th1.1}
Let $0<\alpha<1$ and $u$ be a viscosity solution of
\begin{equation}\label{e-Laplace}
 \left\{ \begin{aligned}
&\Delta u=f ~~&&\mbox{in}~~\Omega\cap B_1;\\
&u=g ~~&&\mbox{on}~~\partial \Omega\cap B_1.
 \end{aligned}\right.
\end{equation}
Suppose that for some integers $k,l\geq 1$, $u\in C^{k,\alpha}(0)$, $f\in C^{k+l-2,\alpha}(0)$, $g\in C^{k+l,\alpha}(0)$ and $\partial \Omega\in C^{l,\alpha}(0)$. Moreover, assume that
\begin{equation*}
u(0)=\cdots=|D^ku(0)|=|Dg(0)|=\cdots=|D^kg(0)|=0.
\end{equation*}

Then $u\in C^{k+l,\alpha}(0)$. That is, there exists $P\in \mathcal{P}_{k+l}$ such that
\begin{equation*}
  \begin{aligned}
&|u(x)-P(x)|\leq C |x|^{k+l+\alpha}\left(\|u\|_{L^{\infty}(\Omega_1)}
+\|f\|_{C^{k+l-2,\alpha}(0)}+\|g\|_{C^{k+l,\alpha}(0)}\right),~\forall ~x\in \Omega\cap B_{1},\\
&|D^{k+1}u(0)|+\cdots+|D^{k+l}u(0)|\leq C\left(\|u\|_{L^{\infty}(\Omega_1)}
+\|f\|_{C^{k+l-2,\alpha}(0)}+\|g\|_{C^{k+l,\alpha}(0)}\right)
  \end{aligned}
\end{equation*}
and
\begin{equation}\label{e1.2}
\begin{aligned}
\Delta P\equiv P_f,~~
\mathbf{\Pi}_{k+l}\left(P(x',P_{\Omega}(x'))\right)\equiv
\mathbf{\Pi}_{k+l}\left(P_g(x',P_{\Omega}(x'))\right),
\end{aligned}
\end{equation}
where $C$ depends only on $n,k,l,\alpha$ and $\|\partial \Omega\cap B_1\|_{C^{l,\alpha}(0)}$.
\end{theorem}

\begin{remark}\label{re1.9}
For the notion of viscosity solutions and related theories, we refer to \cite{MR1351007} and \cite{MR1118699}.
\end{remark}

\begin{remark}\label{re1.22}
Since we consider the pointwise regularity at $0$, $|D^iu(0)|=0$ doesn't imply $|D^ig(0)|=0$ for any $i\geq 1$.
\end{remark}

\begin{remark}\label{re1.12}
In fact, the expression of the polynomial $P$ can be written explicitly:
\begin{equation}\label{e1.5}
  P(x)=P_g(x)+\mathbf{\Pi}_{k+l}\left(\sum_{\mathop{k+1\leq |\sigma|\leq k+l,}\limits_{\sigma_n\geq 1}} \frac{a_{\sigma}}{\sigma !}x^{\sigma-e_n}\left(x_n-P_{\Omega}(x')\right)\right),
\end{equation}
where $a_{\sigma}$ are constants.
\end{remark}

\begin{remark}\label{re1.4}
For equations with coefficients and lower order terms:
\begin{equation*}
  a^{ij}u_{ij}+b^iu_i+cu=f~~\mbox{ in }~~\Omega\cap B_1,
\end{equation*}
we also have $u\in C^{k+l,\alpha}(0)$ ($k\geq 2$) if
\begin{equation*}
a^{ij}\in C^{l-1,\alpha}(0),~b^{i}\in C^{l-2,\alpha}(0)~~\mbox{and}~~c\in C^{l-3,\alpha}(0).
\end{equation*}
\end{remark}

\begin{remark}\label{re1.24}
\Cref{th1.1} can be extended to more general equations, including fully nonlinear uniformly elliptic equations in general forms. For simplicity and clarity, we only consider the Laplace operator in this paper.
\end{remark}

\begin{remark}\label{re1.21}
Note that \Cref{th1.1} can only be stated in the form of pointwise regularity since we can't propose an assumption like $u=|Du|=\cdots=|D^ku|=0$ on $\partial \Omega\cap B_1$.
\end{remark}

\begin{remark}\label{re1.14}
The \cref{e1.2} can be regarded as a polynomial version of \cref{e-Laplace}.
\end{remark}
~\\

As a consequence, we have the following boundary pointwise regularity.
\begin{theorem}\label{co1.1}
Let $0<\alpha<1$ and $u$ be a viscosity solution of
\begin{equation*}
 \left\{ \begin{aligned}
&\Delta u=f ~~&&\mbox{in}~~\Omega\cap B_1;\\
&u=g ~~&&\mbox{on}~~\partial \Omega\cap B_1.
 \end{aligned}\right.
\end{equation*}
Suppose that for some $k\geq 1$, $f\in C^{k-2,\alpha}(0)$, $g\in C^{k,\alpha}(0)$ and $\partial \Omega\in C^{k,\alpha}(0)$.

Then $u\in C^{k,\alpha}(0)$. That is, there exists $P\in \mathcal{P}_k$ such that
\begin{equation*}
  \begin{aligned}
&|u(x)-P(x)|\leq C |x|^{k+\alpha}\left(\|u\|_{L^{\infty}(\Omega_1)}
+\|f\|_{C^{k-2,\alpha}(0)}+\|g\|_{C^{k,\alpha}(0)}\right),~\forall ~x\in \Omega\cap B_{1},\\
&|Du(0)|+\cdots+|D^{k}u(0)|\leq C\left(\|u\|_{L^{\infty}(\Omega_1)}
+\|f\|_{C^{k-2,\alpha}(0)}+\|g\|_{C^{k,\alpha}(0)}\right),
  \end{aligned}
\end{equation*}
where $C$ depends only on $n,k,\alpha$ and $\|\partial \Omega\cap B_1\|_{C^{k,\alpha}(0)}$. Moreover, if $k=1$,
\begin{equation*}
P(x',0)\equiv P_g(x',0);
\end{equation*}
if $k\geq 2$,
\begin{equation*}
\Delta P\equiv P_f,~~
\mathbf{\Pi}_{k}\left(P(x',P_{\Omega}(x'))\right)\equiv
\mathbf{\Pi}_{k}\left(P_g(x',P_{\Omega}(x'))\right),
\end{equation*}

\end{theorem}

\begin{remark}\label{re1.5}
\Cref{co1.1} has been proved in \cite{lian2020pointwise} as a special case. Since the result of \Cref{co1.1} is not well-known and the proof for the Laplace operator is rather simple than that for fully nonlinear equations, we list this result and give a detailed proof in this paper.
\end{remark}

As an application of \Cref{th1.1} to the higher regularity of free boundaries in obstacle-type problems, we have
\begin{theorem}\label{th1.2}
Let $u$ be a viscosity solution of
\begin{equation}\label{e1.1}
\left\{\begin{aligned}
   &\Delta u=1~~~~&&\mbox{in}~~\Omega \cap B_1;\\
   &u=|Du|=0~~~~&&\mbox{on}~~\partial \Omega\cap B_1.\\
\end{aligned}\right.
\end{equation}
Assume that $\partial \Omega\cap B_1\in C^{1,\alpha}$ for some $0<\alpha<1$. Then $u\in C^{\infty} (\bar\Omega\cap B_{1})$ and $\partial \Omega\cap B_1\in C^{\infty}$.
\end{theorem}

\begin{remark}\label{re1.10}
Although we only consider the Poisson equation in this paper, the method is applicable to the higher regularity of free boundaries for fully nonlinear elliptic equations, which have been well studied (see \cite{MR3198649}).
\end{remark}

\begin{remark}\label{re1.11}
For the obstacle-type problem \cref{e1.1}, one can prove the higher regularity starting from that $\partial \Omega$ is Lipschitz continuous with the aid of boundary Harnack inequality (see \cite{MR803243}, \cite[Chapter 6.2]{MR2962060}). Recently, De Silva and Savin \cite{MR4093736} gave an elegant proof of the boundary Harnack inequality for equations in non-divergence form. Hence, we can prove the higher regularity in a simple way from Lipschitz regularity even for fully nonlinear elliptic equations.
\end{remark}

\begin{remark}\label{re1.3}
Usually, the higher regularity of free boundaries is proved in the following way. First, by a proper partial hodograph-Legendre transformation (see \cite{MR0440187}), \cref{e1.1} is transformed to an elliptic equation on a flat boundary and the free boundary $\partial \Omega\cap B_1$ is represented by some function relation. Note that even for the Laplace operator, the transformed equation is a fully nonlinear elliptic equation in general. Hence, in the second step, one need to apply the theory for fully nonlinear elliptic equations (see \cite{MR0125307}) to obtain the higher regularity of solutions and free boundaries.

In 2015, De Silva and Savin \cite{MR3393271} gave a direct and simple proof of higher regularity of solutions and free boundaries based on a higher order boundary Harnack inequality. However, this method is not applicable to the fully nonlinear elliptic equations.
\end{remark}
~\\

With respect to the boundary pointwise regularity for oblique derivative problems, we have
\begin{theorem}\label{th1.3-1}
Let $0<\alpha<1$ and $u$ be a viscosity solution of
\begin{equation*}
 \left\{ \begin{aligned}
&\Delta u=f ~~&&\mbox{in}~~\Omega\cap B_1;\\
&\beta\cdot Du=g ~~&&\mbox{on}~~\partial \Omega\cap B_1,
 \end{aligned}\right.
\end{equation*}
where $\beta$ is a vector valued function and satisfies the obliqueness condition: for some positive constant $a_0$,
\begin{equation}\label{e1.4}
\beta_n=\beta\cdot e_n\geq a_0.
\end{equation}
Suppose that $f\in C^{-1,\alpha}(0)$, $g\in C^{\alpha}(0)$, $\beta\in C^{\alpha}(0)$ and $[\partial \Omega\cap B_1]_{C^{0,1}(0)}\leq \delta$, where $\delta>0$ depending only on $n,\alpha$ and $a_0$.

Then $u\in C^{1,\alpha}(0)$. That is, there exists $P\in \mathcal{P}_1$ such that
\begin{equation*}
  \begin{aligned}
&|u(x)-P(x)|\leq C |x|^{1+\alpha}\left(\|u\|_{L^{\infty}(\Omega_1)}
+\|f\|_{C^{-1,\alpha}(0)}+\|g\|_{C^{\alpha}(0)}\right),~~\forall ~x\in \Omega\cap B_{1},\\
&|Du(0)|\leq C\left(\|u\|_{L^{\infty}(\Omega_1)}
+\|f\|_{C^{-1,\alpha}(0)}+\|g\|_{C^{\alpha}(0)}\right)
  \end{aligned}
\end{equation*}
and
\begin{equation*}
  \beta(0)\cdot Du(0)=g(0).
\end{equation*}
where $C$ depends only on $n,\alpha$ and $a_0$.
\end{theorem}

\begin{remark}\label{re1.20}
We refer the reader to \cite{MR3780142} and \cite{MR2254613} for the notion of viscosity solutions for oblique derivative problems.
\end{remark}

\begin{remark}\label{re1.6}
Without loss of generality, for oblique derivative problems, we always assume that $\|\beta\|_{L^{\infty}(\partial \Omega\cap B_1)}\leq 1$ and \cref{e1.4} holds with the fixed constant $a_0$ throughout this paper.
\end{remark}

\begin{remark}\label{re1.7}
\Cref{th1.3-1} shows that the boundary $C^{1,\alpha}$ regularity holds on Lipschitz domains with a small Lipschitz constant.
\end{remark}
~\\

For higher regularity, we have
\begin{theorem}\label{th1.3}
Let $0<\alpha<1$ and $u$ be a viscosity solution of
\begin{equation*}
 \left\{ \begin{aligned}
&\Delta u=f ~~&&\mbox{in}~~\Omega\cap B_1;\\
&\beta\cdot Du=g ~~&&\mbox{on}~~\partial \Omega\cap B_1.
 \end{aligned}\right.
\end{equation*}
Suppose that for some integers $k,l\geq 1$, $u\in C^{k,\alpha}(0)$, $f\in C^{k+l-2,\alpha}(0)$, $g\in C^{k+l-1,\alpha}(0)$, $\beta\in C^{l-1,\alpha}(0)$ and $\partial \Omega\in C^{l,\alpha}(0)$. Moreover, assume that
\begin{equation*}
u(0)=\cdots=|D^ku(0)|=g(0)=\cdots=|D^{k-1}g(0)|=0.
\end{equation*}

Then $u\in C^{k+l,\alpha}(0)$. That is, there exists $P\in \mathcal{P}_{k+l}$
such that
\begin{equation*}
  \begin{aligned}
&|u(x)-P(x)|\leq C |x|^{k+l+\alpha}\left(\|u\|_{L^{\infty}(\Omega_1)}
+\|f\|_{C^{k+l-2,\alpha}(0)}+\|g\|_{C^{k+l-1,\alpha}(0)}\right),~~\forall ~x\in \Omega\cap B_{1},\\
&|D^{k+1}u(0)|+\cdots+|D^{k+l}u(0)|\leq C\left(\|u\|_{L^{\infty}(\Omega_1)}
+\|f\|_{C^{k+l-2,\alpha}(0)}+\|g\|_{C^{k+l-1,\alpha}(0)}\right)
  \end{aligned}
\end{equation*}
and
\begin{equation}\label{e1.3}
\Delta P\equiv P_f,~~\mathbf{\Pi}_{k+l-1}\bigg((P_{\beta}\cdot DP)(x',P_{\Omega}(x'))\bigg)\equiv \mathbf{\Pi}_{k+l-1}\bigg(P_g(x',P_{\Omega}(x'))\bigg),
\end{equation}
where $C$ depends only on $n,k,l,\alpha,a_0$, $\|\beta\|_{C^{l-1,\alpha}(0)}$ and $\|\partial \Omega\cap B_1\|_{C^{l,\alpha}(0)}$.
\end{theorem}

Similar to the Dirichlet problem, we have the following higher order boundary pointwise regularity.
\begin{theorem}\label{th1.5}
Let $0<\alpha<1$ and $u$ be a viscosity solution of
\begin{equation*}
 \left\{ \begin{aligned}
&\Delta u=f ~~&&\mbox{in}~~\Omega\cap B_1;\\
&\beta\cdot Du=g ~~&&\mbox{on}~~\partial \Omega\cap B_1.
 \end{aligned}\right.
\end{equation*}
Suppose that for some $k\geq 2$, $f\in C^{k-2,\alpha}(0)$, $g\in C^{k-1,\alpha}(0)$, $\beta\in C^{k-1,\alpha}(0)$ and $\partial \Omega\in C^{k-1,\alpha}(0)$.

Then $u\in C^{k,\alpha}(0)$. That is, there exists $P\in \mathcal{P}_k$ such that
\begin{equation*}
  \begin{aligned}
&|u(x)-P(x)|\leq C |x|^{k+\alpha}\left(\|u\|_{L^{\infty}(\Omega_1)}
+\|f\|_{C^{k-2,\alpha}(0)}+\|g\|_{C^{k-1,\alpha}(0)}\right),~~\forall ~x\in \Omega\cap B_{1}\\
&|Du(0)|+\cdots+|D^{k}u(0)|\leq C\left(\|u\|_{L^{\infty}(\Omega_1)}
+\|f\|_{C^{k-2,\alpha}(0)}+\|g\|_{C^{k-1,\alpha}(0)}\right)
  \end{aligned}
\end{equation*}
and
\begin{equation*}
\Delta P\equiv P_f,~~\mathbf{\Pi}_{k-1}\bigg((P_{\beta}\cdot DP)(x',P_{\Omega}(x'))\bigg)\equiv \mathbf{\Pi}_{k-1}\bigg(P_g(x',P_{\Omega}(x'))\bigg),
\end{equation*}
where $C$ depends only on $n,k,\alpha,a_0$, $\|\beta\|_{C^{k-1,\alpha}(0)}$ and $\|\partial \Omega\cap B_1\|_{C^{k-1,\alpha}(0)}$.
\end{theorem}

\begin{remark}\label{re1.23}
Similar to the Dirichlet problems, one can prove corresponding pointwise boundary regularity for fully nonlinear elliptic equations.
\end{remark}

As an application of the above boundary pointwise regularity to the regularity of free boundaries in one phase problems, we have
\begin{theorem}\label{th1.4}
Let $0<\alpha<1$ and $u$ be a viscosity solution of
\begin{equation*}
\left\{\begin{aligned}
&\Delta u=1~~~~&&\mbox{in}~~\Omega \cap B_1;\\
&u=0~~~~&&\mbox{on}~~\partial \Omega\cap B_1;\\
&|Du|=1~~~~&&\mbox{on}~~\partial \Omega\cap B_1.
\end{aligned}\right.
\end{equation*}
Assume that $\partial \Omega\cap B_1\in C^{1,\alpha}$. Then $u\in C^{\infty} (\bar\Omega\cap B_{1})$ and $\partial \Omega\cap B_1\in C^{\infty}$.
\end{theorem}

In \Cref{S2}, we prove the boundary pointwise regularity for Dirichlet problems and the higher regularity for free boundaries of obstacle-type problems. \Cref{S3} is devoted to the oblique derivative problems and higher regularity for free boundaries in one phase problems. The results of this paper show the underlying relation between the regularity of solutions and the regularity of boundaries. The proofs demonstrate how overdetermined conditions lead to higher regularity. Notations used in this paper are listed below, most of which are standard.
\begin{notation}\label{no1.1}
\begin{enumerate}~~\\
\item $\{e_i\}^{n}_{i=1}$: the standard basis of $\mathbb{R}^n$, i.e., $e_i=(0,...0,\underset{i^{th}}{1},0,...0)$.
\item $x'=(x_1,x_2,...,x_{n-1})$ and $x=(x_1,...,x_n)=(x',x_n)$ .
\item $|x|:=\left(\sum_{i=1}^{n} x_i^2\right)^{1/2}$ for $x\in \mathbb{R}^n$.
\item $\mathbb{R}^n_+=\{x\in \mathbb{R}^n\big|x_n>0\}$.
\item $B_r(x_0)=B(x_0,r)=\{x\in \mathbb{R}^{n}\big| |x-x_0|<r\}$, $B_r=B_r(0)$, $B_r^+(x_0)=B_r(x_0)\cap \mathbb{R}^n_+$ and $B_r^+=B^+_r(0)$.
\item $T_r(x_0)\ =\{(x',0)\in \mathbb{R}^{n}\big| |x'-x_0'|<r\}$ \mbox{ and } $T_r=T_r(0)$.
\item $A^c$: the complement of $A$; $\bar A $: the closure of $A$, where $ A\subset \mathbb{R}^n$.
\item $\mathrm{diam}(A)$: the diameter of $A$ and $\mathrm{dist}(A,B)$: the distance between $A$ and $B$, where $ A,B\subset \mathbb{R}^n$.
\item $\Omega_r=\Omega\cap B_r$ and $(\partial\Omega)_r=\partial\Omega\cap B_r$.
\item $\varphi _i=D_i \varphi=\partial \varphi/\partial x _{i}$, $\varphi _{ij}=D_{ij}\varphi =\partial ^{2}\varphi/\partial x_{i}\partial x_{j}$ and we also use similar notations for higher order derivatives.
\item $D^0\varphi =\varphi$, $D \varphi=(\varphi_1 ,...,\varphi_{n} )$ and $D^2 \varphi =\left(\varphi _{ij}\right)_{n\times n}$ etc.
\item $|D^k\varphi |= \left(\sum_{|\sigma|=k}|D^{\sigma}\varphi| ^2\right)^{1/2}$ for $k\geq 1$, where the standard multi-index notation is used.
\item $\mathcal{P}_k (k\geq 0):$ the set of polynomials of degree less than or equal to $k$. That is, any $P\in \mathcal{P}_k$ can be written as
\begin{equation*}
P(x)=\sum_{|\sigma|\leq k}\frac{a_{\sigma}}{\sigma!}x^{\sigma}
\end{equation*}
where $a_{\sigma}$ are constants. Define
\begin{equation*}
\|P\|= \sum_{|\sigma|\leq k}|a_{\sigma}|.
\end{equation*}
\item $\mathcal{HP}_k (k\geq 0):$ the set of homogeneous polynomials of degree $k$. That is, any $P\in \mathcal{HP}_k$ can be written as
\begin{equation*}
P(x)=\sum_{|\sigma|= k}\frac{a_{\sigma}}{\sigma!}x^{\sigma}.
\end{equation*}
\item $\mathbf{\Pi}_k:$ The projection from $\mathcal{P}_l$ to $\mathcal{P}_k$ for $l\geq k$. That is, if $P\in \mathcal{P}_l$ is written as
\begin{equation*}
P(x)=\sum_{|\sigma|\leq l}\frac{a_{\sigma}}{\sigma!}x^{\sigma},
\end{equation*}
then
\begin{equation*}
\mathbf{\Pi}_kP(x)=\sum_{|\sigma|\leq k}\frac{a_{\sigma}}{\sigma!}x^{\sigma}.
\end{equation*}
\end{enumerate}
\end{notation}
~\\

\section{Dirichlet problem and application to the obstacle problem}\label{S2}
In this section, we give the proofs of \Cref{th1.1}, \Cref{co1.1} and \Cref{th1.2}. We start with the following result (see \cite[Corollary 4.3]{MR4088470}):
\begin{lemma}\label{l-2.1}
Let $0<\alpha<1$ and $u$ be a viscosity solution of
\begin{equation*}
\left\{\begin{aligned}
&\Delta u=f&& ~~\mbox{in}~~\Omega_1;\\
&u=g&& ~~\mbox{on}~~(\partial \Omega)_1.
\end{aligned}\right.
\end{equation*}
Suppose that $u\in C^{1,\alpha}(0)$, $f\in C^{\alpha}(0)$, $g\in C^{2,\alpha}(0)$ and $(\partial \Omega)_1\in C^{1,\alpha}(0)$. Moreover, assume that
\begin{equation*}
u(0)=|Du(0)|=|Dg(0)|=0.
\end{equation*}

Then $u\in C^{2,\alpha}(0)$. That is, there exists $P\in \mathcal{HP}_2$ such that
\begin{equation*}
  \begin{aligned}
&|u(x)-P(x)|\leq C |x|^{2+\alpha}\left(\|u\|_{L^{\infty}(\Omega_{1})}+\|f\|_{C^{\alpha}(0)}
+\|g\|_{C^{2,\alpha}(0)}\right),~~ \forall ~x\in \Omega_{1},\\
&|D^{2}u(0)|\leq C\left(\|u\|_{L^{\infty}(\Omega_{1})}+\|f\|_{C^{\alpha}(0)}
+\|g\|_{C^{2,\alpha}(0)}\right)
  \end{aligned}
\end{equation*}
and
\begin{equation}\label{e2.11}
  \begin{aligned}
    \Delta P&=f(0),~~P(x',0)\equiv P_g(x',0),
  \end{aligned}
\end{equation}
where $C$ depends only on $n,\alpha$ and $\|(\partial \Omega)_1\|_{C^{1,\alpha}(0)}$.
\end{lemma}

Next, we prove a generalized version of \Cref{l-2.1}:
\begin{lemma}\label{l-2.2}
Let $0<\alpha<1$ and $u$ be a viscosity solution of
\begin{equation*}
\left\{\begin{aligned}
&\Delta u=f&& ~~\mbox{in}~~\Omega_1;\\
&u=g&& ~~\mbox{on}~~(\partial \Omega)_1.
\end{aligned}\right.
\end{equation*}
Suppose that $u\in C^{k,\alpha}(0) (k\geq 1)$, $f\in C^{k-1,\alpha}(0)$, $g\in C^{k+1,\alpha}(0)$ and $(\partial \Omega)_1\in C^{1,\alpha}(0)$. Moreover, assume that
\begin{equation*}
u(0)=\cdots=|D^ku(0)|=|Dg(0)|\cdots=|D^kg(0)|=0.
\end{equation*}

Then $u\in C^{k+1,\alpha}(0)$. That is, there exists $P\in \mathcal{HP}_{k+1}$ such that
\begin{equation*}
  \begin{aligned}
&|u(x)-P(x)|\leq C |x|^{k+1+\alpha}\left(\|u\|_{L^{\infty}(\Omega_{1})}+\|f\|_{C^{k-1,\alpha}(0)}
+\|g\|_{C^{k+1,\alpha}(0)}\right),~~\forall ~x\in \Omega_{1},\\
&|D^{k+1}u(0)|\leq C\left(\|u\|_{L^{\infty}(\Omega_{1})}+\|f\|_{C^{k-1,\alpha}(0)}
+\|g\|_{C^{k+1,\alpha}(0)}\right)
  \end{aligned}
\end{equation*}
and
\begin{equation}\label{e2.9}
  \begin{aligned}
    \Delta P&\equiv P_f,~~P(x',0)\equiv P_g(x',0),
  \end{aligned}
\end{equation}
where $C$ depends only on $n,k,\alpha$ and $\|(\partial \Omega)_1\|_{C^{1,\alpha}(0)}$.
\end{lemma}

In the following proof, we will use a kind of homogeneous polynomial in a special form. We call $Q\in \mathcal{HP}_k$ a $k$-form ($k\geq 1$) if $Q$ can be written as
\begin{equation*}
Q(x)=\sum_{|\sigma|=k,\sigma_n\geq 1}\frac{a_{\sigma}}{\sigma !}x^{\sigma}.
\end{equation*}
That is, $x_n$ appears in the expression of $Q$ at least one time (thus $Q\equiv 0$ on $T_1$), which turns out to be vital for the boundary regularity. In fact, $P(x',0)\equiv P_g(x',0)$ in \cref{e2.9} indicates that $P(x)-P_g(x)$ is a $(k+1)$-form.

We prove \Cref{l-2.2} by induction. For $k=1$, the lemma reduces to \Cref{l-2.1}. Suppose that the lemma holds for $k\leq k_0-1$ and we need to prove the lemma for $k=k_0$. First, we prove a key step towards the conclusion of \Cref{l-2.2}.

\begin{lemma}\label{l-2.3}
Let $1\leq k\leq k_0$, $0<\alpha<1$ and $u\in C^{k,\alpha}(0)$ be a viscosity solution of
\begin{equation*}
\left\{\begin{aligned}
&\Delta u+P=f&& ~~\mbox{in}~~\Omega_1;\\
&u=g&& ~~\mbox{on}~~(\partial \Omega)_1,
\end{aligned}\right.
\end{equation*}
where $P\in \mathcal{HP}_{k-1}$. Suppose that
\begin{equation*}
  \begin{aligned}
    &\|u\|_{L^{\infty}(\Omega_1)}\leq 1, u(0)=\cdots=|D^ku(0)|=0,\\
    &|f(x)|\leq \delta|x|^{k-2+\alpha}, ~~\forall x\in \Omega_1,\\
    &|g(x)|\leq \delta|x|^{k+\alpha},~~\forall x\in (\partial \Omega)_1,\\
    &\|(\partial \Omega)_1\|_{C^{1,\alpha}(0)} \leq \delta,\\
    &\|P\|\leq 1,
  \end{aligned}
\end{equation*}
where $\delta>0$ depends only on $n,k$ and $\alpha$.

Then there exists a $(k+1)$-form $Q$ such that
\begin{equation*}
  \begin{aligned}
    & \|u-Q\|_{L^{\infty}(\Omega_{\eta})}\leq \eta^{k+1+\alpha},\\
    &\|Q\|\leq C_0,\\
    &\Delta Q+P\equiv 0,\\
  \end{aligned}
\end{equation*}
where $C_0$ depends only on $n$ and $k$, and $\eta$ depends also on $\alpha$.
\end{lemma}
\proof We prove the lemma by contradiction. Suppose that the conclusion is false. Then there exist $0<\alpha<1$ and sequences of $u_m,f_m,g_m,\Omega_m,P_m$ ($m\geq 1$) satisfying $u_m\in C^{k,\alpha}(0)$ and
\begin{equation*}
\left\{\begin{aligned}
&\Delta u_m+P_m=f_m&& ~~\mbox{in}~~\Omega_m\cap B_1;\\
&u_m=g_m&& ~~\mbox{on}~~\partial \Omega_m\cap B_1.
\end{aligned}\right.
\end{equation*}
In addition,
\begin{equation*}
  \begin{aligned}
    &\|u_m\|_{L^{\infty}(\Omega_m\cap B_1)}\leq 1, u_m(0)=\cdots=|D^ku_m(0)|=0,\\
    &|f_m(x)|\leq |x|^{k-2+\alpha}/m, ~~\forall x\in \Omega_m\cap B_1,\\
    &|g_m(x)|\leq |x|^{k+\alpha}/m,~~\forall x\in \partial \Omega_m\cap B_1,\\
    &\|\partial \Omega_m\cap B_1\|_{C^{1,\alpha}(0)} \leq 1/m,\\
    &\|P_m\|\leq 1.
  \end{aligned}
\end{equation*}
But for any $(k+1)$-form $Q$ satisfying $\|Q\|\leq C_0$ and $\Delta Q+P_m\equiv 0$, we have
\begin{equation}\label{e.lCka.1-mu}
  \|u_m-Q\|_{L^{\infty}(\Omega_{m}\cap B_{\eta})}> \eta^{k+1+\alpha},
\end{equation}
where $C_0$ is to be specified later and $0<\eta<1$ is taken small such that
\begin{equation}\label{e.lCka.2-mu}
C_0\eta^{1-\alpha}<1/2.
\end{equation}

Clearly, $u_m$ are uniformly bounded ($\|u_m\|_{L^{\infty}(\Omega_m\cap B_1)}\leq 1$). Moreover, $u_m$ are equicontinuous (see \cite[Lemma 2.7]{MR4088470}). Hence, there exist $\tilde u:B_1^+\cup T_1\rightarrow \mathbb{R}$ and $\tilde P\in \mathcal{HP}_{k-1}$ such that $u_m\rightarrow\tilde u$ uniformly in compact subsets of $ B_1^+\cup T_1$, $P_m\rightarrow \tilde P$ and
\begin{equation*}
\left\{\begin{aligned}
&\Delta \tilde u+\tilde P=0&& ~~\mbox{in}~~B_{1}^+;\\
&\tilde u=0&& ~~\mbox{on}~~T_{1}.
\end{aligned}\right.
\end{equation*}

By the boundary $C^{k,\alpha}$ estimate for $u_m$ (\Cref{l-2.2} for $k-1$ since $k\leq k_0$) and noting $u_m(0)=\cdots=|D^{k}u_m(0)|=0$, we have
\begin{equation*}
\|u_m\|_{L^{\infty }(\Omega_m\cap B_r)}\leq Cr^{k+\alpha} ,~~~~\forall~0<r<1.
\end{equation*}
Since $u_m$ converges to $u$ uniformly,
\begin{equation*}
\|\tilde u\|_{L^{\infty }(B_r^+)}\leq Cr^{k+\alpha}, ~~~~\forall~0<r<1.
\end{equation*}
Hence, $\tilde u(0)=\cdots=|D^{k}\tilde u(0)|=0$. By the boundary estimate for $\tilde u$ on a flat boundary, there exists a $(k+1)$-form $\tilde{Q}$ such that
\begin{equation}\label{e.cka-5}
\begin{aligned}
  &|\tilde u(x)-\tilde{Q}(x)|\leq C_0 |x|^{k+2}, ~~\forall ~x\in B_{1}^+,\\
  &\Delta \tilde{Q}+\tilde P\equiv 0,\\
  &\|\tilde{Q}\|\leq C_0/2,
\end{aligned}
\end{equation}
where $C_0$ depends only on $n$ and $k$.

Since $P_m\rightarrow \tilde{P}$, there exist $(k+1)$-forms $\tilde{Q}_m$ such that $\Delta (\tilde{Q}+ \tilde{Q}_m)+P_m\equiv 0$ and $\|\tilde{Q}_m\|\rightarrow 0$ as $m\rightarrow \infty$. Thus, \cref{e.lCka.1-mu} holds for $Q=\tilde{Q}+\tilde{Q}_m$. Let $m\rightarrow \infty$ in \cref{e.lCka.1-mu} and we have
\begin{equation*}
    \|\tilde u-\tilde{Q}\|_{L^{\infty}(B_{\eta}^+)}\geq \eta^{k+1+\alpha},
\end{equation*}
However, by \cref{e.lCka.2-mu} and \cref{e.cka-5},
\begin{equation*}
  \|\tilde u-\tilde{Q}\|_{L^{\infty}(B_{\eta}^+)}\leq \eta^{k+1+\alpha}/2,
\end{equation*}
which is a contradiction.  ~\qed~\\

Now, we give the\\
\noindent\textbf{Proof of \Cref{l-2.2}.} Since we have assumed that \Cref{l-2.2} holds for $k_0-1$, $u\in C^{k_0,\alpha}(0)$. By induction, we only need to prove \Cref{l-2.2} for $k_0$, i.e., $u\in C^{k_0+1,\alpha}(0)$. Without loss of generality, by a proper transformation, we assume that
\begin{equation}\label{e2.7}
\left\{\begin{aligned}
&\Delta u+P=f&& ~~\mbox{in}~~\Omega_1;\\
&u=g&& ~~\mbox{on}~~(\partial \Omega)_1
\end{aligned}\right.
\end{equation}
for some $P\in \mathcal{HP}_{k_0-1}$ and
\begin{equation}\label{e2.1}
  \begin{aligned}
    &\|u\|_{L^{\infty}(\Omega_1)}\leq 1, u(0)=\cdots=|D^{k_0}u(0)|=0,\\
    &|f(x)|\leq \delta |x|^{k_0-1+\alpha}, ~~\forall x\in \Omega_1,\\
    &|g(x)|\leq \delta |x|^{k_0+1+\alpha}/2,~~\forall x\in (\partial \Omega)_1,\\
    &\|(\partial \Omega)_1\|_{C^{1,\alpha}(0)} \leq \delta/(2C_1),\\
    &\|P\|\leq 1,
  \end{aligned}
\end{equation}
where $\delta$ is as in \Cref{l-2.3} (with $k=k_0$) and $C_1$ depending only on $n,k_0$ and $\alpha$ is to be chosen later.

Indeed, let $u_1=u/K$ where
\begin{equation*}
K=2(\|u\|_{L^{\infty}(\Omega_1)}+\|f\|_{C^{k_0-1,\alpha}(0)}
+\|g\|_{C^{k_0+1,\alpha}(0)}).
\end{equation*}
Then $u_1$ satisfies
\begin{equation*}
\left\{\begin{aligned}
&\Delta u_1=f_1&& ~~\mbox{in}~~\Omega_1;\\
&u_1=g_1&& ~~\mbox{on}~~(\partial \Omega)_1,
\end{aligned}\right.
\end{equation*}
where $f_1=f/K$ and $g_1=g/K$.

Next, let $u_2=u_1-P_{g_1}$. Since $g(0)=|Dg(0)|\cdots=|D^{k_0}g(0)|=0$, $P_g\in \mathcal{HP}_{k_0+1}$. In addition, by $u(0)=|Du(0)|\cdots=|D^{k_0}u(0)|=0$, we have $f(0)=\cdots=|D^{k_0-2}f(0)|=0$ (see \cref{e2.9}). Hence, $P_f\in \mathcal{HP}_{k_0-1}$. Then $u_2$ satisfies
\begin{equation*}
\left\{\begin{aligned}
&\Delta u_2+P=f_2&& ~~\mbox{in}~~\Omega_1;\\
&u_2=g_2&& ~~\mbox{on}~~(\partial \Omega)_1,
\end{aligned}\right.
\end{equation*}
where
\begin{equation*}
  |f_2(x)|=|f_1(x)-P_{f_1}(x)|\leq |x|^{k_0-1+\alpha}, ~~\forall x\in \Omega_1,
\end{equation*}
\begin{equation*}
  |g_2(x)|=|g_1(x)-P_{g_1}(x)|\leq |x|^{k_0+1+\alpha}, ~~\forall x\in (\partial\Omega)_1.
\end{equation*}
In addition, $P=\Delta P_{g_1}-P_{f_1}\in\mathcal{HP}_{k_0-1}$ and
\begin{equation*}
  \|P\|\leq C,
\end{equation*}
where $C$ depends only on $n$ and $k_0$.

Finally, let $y=x/\rho$ for $\rho>0$ and $\tilde{u}(y)=u_2(x)$. Then $\tilde{u}$ satisfies
\begin{equation*}
\left\{\begin{aligned}
&\Delta \tilde{u}+\tilde{P}=\tilde f&& ~~\mbox{in}~~\tilde\Omega_1;\\
&\tilde{u}=\tilde{g}&& ~~\mbox{on}~~(\partial \tilde\Omega)_1,
\end{aligned}\right.
\end{equation*}
where
\begin{equation*}
\tilde{f}(y)=\rho^2f_2(x),~\tilde{g}(y)=g_2(x),~\tilde{P}(y)=\rho^{k+1}P(x),~
\tilde{\Omega}=\Omega/\rho.
\end{equation*}
Hence,
\begin{equation*}
  \begin{aligned}
    &\|\tilde u\|_{L^{\infty}(\tilde\Omega_1)}=\|u_2\|_{L^{\infty}(\Omega_{\rho})}
    \leq \|u_1\|_{L^{\infty}(\Omega_1)}+\|P_{g_1}\|\leq 1, \tilde{u}(0)=\cdots=|D^{k_0}{u}(0)|=0,\\
    &|\tilde{f}(y)|\leq \rho^{k_0+1+\alpha} |y|^{k_0-1+\alpha}, ~~\forall y\in \tilde\Omega_1,\\
    &|\tilde{g}(y)|\leq \rho^{k_0+1+\alpha} |y|^{k_0+1+\alpha},~~\forall y\in (\partial \tilde\Omega)_1,\\
    &\|(\partial \tilde\Omega)_1\|_{C^{1,\alpha}(0)}\leq \rho^{\alpha}\|(\partial\Omega)_1\|_{C^{1,\alpha}(0)},\\
    &\|\tilde{P}\|= \rho^{k_0+1+\alpha}\|P\|\leq  \rho^{k_0+1+\alpha}C.
  \end{aligned}
\end{equation*}

Therefore, by taking $\rho$ small enough (depending only on $n,k_0,\alpha$ and $\|(\partial\Omega)_1\|_{C^{1,\alpha}(0)}$), the assumptions \cref{e2.7} and \cref{e2.1} for $\tilde{u}$ can be guaranteed. Then the regularity of $u$ can be derived from that of $\tilde{u}$. Hence, without loss of generality, we assume that \cref{e2.7} and \cref{e2.1} hold for $u$.

To prove \Cref{l-2.2} for $k_0$, we only need to show that there exist a sequence of $(k_0+1)$-forms $Q_m$ ($m\geq 0$) such that for all $m\geq 1$,
\begin{equation}\label{e2.2}
\|u-Q_m\|_{L^{\infty }(\Omega _{\eta^{m}})}\leq \eta ^{m(k_0+1+\alpha )},
\end{equation}
\begin{equation}\label{e2.3}
\Delta Q_m+P\equiv 0
\end{equation}
and
\begin{equation}\label{e2.5}
\|Q_m-Q_{m-1}\|\leq C_0\eta ^{m\alpha},
\end{equation}
where $C_0$ and $\eta$ are the constants as in \Cref{l-2.3}.

We prove the above by induction. For $m=1$, by \Cref{l-2.3} and setting $Q_0\equiv 0$, the conclusion holds clearly. Suppose that the conclusion holds for $m$. We need to prove that the conclusion holds for $m+1$.

Let $r=\eta^{m}$, $y=x/r$ and
\begin{equation}\label{e2.6}
  \tilde{u}(y)=\frac{u(x)-Q_m(x)}{r^{k_0+1+\alpha}}.
\end{equation}
Then $\tilde{u}$ satisfies
\begin{equation*}
\left\{\begin{aligned}
&\Delta \tilde{u}=\tilde{f}&& ~~\mbox{in}~~\tilde{\Omega}\cap B_1;\\
&\tilde{u}=\tilde{g}&& ~~\mbox{on}~~\partial \tilde{\Omega}\cap B_1,
\end{aligned}\right.
\end{equation*}
where
\begin{equation*}
 \tilde{f}(y)=\frac{f(x)}{r^{k_0-1+\alpha}},~~
 \tilde{g}(y)=\frac{g(x)-Q_m(x)}{r^{k_0+1+\alpha}},~~
 \tilde{\Omega}=\frac{\Omega}{r}.
\end{equation*}

From \cref{e2.5}, there exists $C_1$ depends only on $n,k_0$ and $\alpha$ such that $\|Q_i\|\leq C_1$
($\forall~0\leq i\leq m$). By combining that $Q_m$ is a $(k_0+1)$-form and the definition of $\partial \Omega\in C^{1,\alpha}(0)$ (see \cref{e-re} and \cref{e-re2}), we have
\begin{equation}\label{e2.8}
\begin{aligned}
|Q_m(x)|&\leq \|Q_m\||x|^{k_0}|x_n|\\
&\leq C_1|x|^{k_0}\|(\partial \Omega)_1\|_{C^{1,\alpha}(0)}|x'|^{1+\alpha}\\
&\leq C_1\|(\partial \Omega)_1\|_{C^{1,\alpha}(0)}|x|^{k_0+1+\alpha},~\forall ~x\in (\partial \Omega)_1.
\end{aligned}
\end{equation}
Therefore,
\begin{equation*}
  \begin{aligned}
& \|\tilde{u}\|_{L^{\infty}(\tilde{\Omega}\cap B_1)}\leq 1, ~(\mathrm{by}~ \cref{e2.2} ~\mbox{and}~ \cref{e2.6})\\
& |\tilde{f}(y)|=\frac{|f(x)|}{r^{k_0-1+\alpha}}\leq \delta|y|^{k_0-1+\alpha},~~\forall y\in \tilde\Omega_1, ~(\mathrm{by}~ \cref{e2.1})\\
&|\tilde{g}(y)|\leq \frac{1}{r^{k_0+1+\alpha}}\left(|g(x)|+|Q_m(x)|\right)\\
&~~~~~~~\leq \frac{1}{r^{k_0+1+\alpha}}\left(\frac{\delta}{2}|x|^{k_0+1+\alpha}
      +C_1\cdot\frac{\delta}{2C_1}|x|^{k_0+1+\alpha}\right)\\
&~~~~~~~\leq \delta |y|^{k_0+1+\alpha},~~\forall y\in (\partial \tilde\Omega)_1, ~(\mathrm{by}~ \cref{e2.1} \mbox{ and } \cref{e2.8})\\
&\|\partial \tilde{\Omega}\cap B_1\|_{C^{1,\alpha}(0)} \leq \delta r^{\alpha} \leq \delta.   ~(\mathrm{by}~ \cref{e2.1})\\
  \end{aligned}
\end{equation*}

By virtue of \Cref{l-2.3}, there exists a $(k_0+1)$-form $\tilde{Q}$ such that
\begin{equation*}
\begin{aligned}
    &\|\tilde{u}-\tilde{Q}\|_{L^{\infty }(\tilde{\Omega} _{\eta})}\leq \eta ^{k_0+1+\alpha},\\
    &\Delta\tilde{Q}\equiv 0,\\
    &\|\tilde{Q}\|\leq C_0.
\end{aligned}
\end{equation*}
Let $Q_{m+1}(x)=Q_m(x)+r^{k_0+1+\alpha}\tilde{Q}(y)=Q_m(x)+r^{\alpha}\tilde{Q}(x)$. Then \cref{e2.3} and \cref{e2.5} hold for $m+1$. Recalling \cref{e2.6}, we have
\begin{equation*}
  \begin{aligned}
&\|u-Q_{m+1}\|_{L^{\infty}(\Omega_{\eta^{m+1}})}\\
&= \|u-Q_m-r^{\alpha}\tilde{Q}\|_{L^{\infty}(\Omega_{\eta r})}\\
&= \|r^{k_0+1+\alpha}\tilde{u}-r^{k_0+1+\alpha}\tilde{Q}\|_{L^{\infty}(\tilde{\Omega}_{\eta})}\\
&\leq r^{k_0+1+\alpha}\eta^{k_0+1+\alpha}\\
&=\eta^{(m+1)(k_0+1+\alpha)}.
  \end{aligned}
\end{equation*}
Hence, \cref{e2.2} holds for $m+1$. By induction, the proof is completed.\qed~\\

\begin{remark}\label{re2.1}
By checking the proof, the condition on $\partial \Omega$ is exactly used for estimating $Q_m$ on $\partial \Omega$ (see \cref{e2.8}). Hence, if the derivatives of $u$ vanish, $Q_m$ will be a higher order homogenous polynomial. This leads to a lower regularity assumption on $\partial \Omega$. This is why we can obtain the $C^{k_0+1,\alpha}$ regularity based only on $\partial \Omega\in C^{1,\alpha}$.
\end{remark}
~\\

Now, we can prove \Cref{th1.1} based on \Cref{l-2.2}.\\
\noindent\textbf{Proof of \Cref{th1.1}.} Throughout this proof, $C$ always denotes a constant depending only on $n,k,l,\alpha$ and $\|(\partial \Omega)_1\|_{C^{l,\alpha}(0)}$. Without loss of generality, we assume
\begin{equation*}
\|u\|_{L^{\infty}(\Omega_1)}
+\|f\|_{C^{k+l-2,\alpha}(0)}+\|g\|_{C^{k+l,\alpha}(0)}\leq 1.
\end{equation*}
Since $g\in C^{k+l,\alpha}(0)$,
\begin{equation*}
  |g(x)-P_g(x)|\leq |x|^{k+l+\alpha},~~~\forall ~x\in (\partial \Omega)_1.
\end{equation*}
Set $u_1=u-P_g$ and then $u_1$ is a viscosity solution of
\begin{equation*}
\left\{\begin{aligned}
&\Delta u_1=f_1&& ~~\mbox{in}~~\Omega_1;\\
&u_1=g_1&& ~~\mbox{on}~~(\partial \Omega)_1,
\end{aligned}\right.
\end{equation*}
where $f_1=f-\Delta P_g$ and $g_1=g-P_g$. Hence,
\begin{equation}\label{e2.10}
g_1(0)=|Dg_1(0)|=\cdots=|D^{k+l}g_1(0)|=0.
\end{equation}

By \Cref{l-2.2}, $u_1\in  C^{k+1,\alpha}(0)$ and there exists a $(k+1)$-form $Q_{k+1}$ such that
\begin{equation*}
|u_1(x)-Q_{k+1}(x)|\leq C |x|^{k+1+\alpha},~~\forall ~x\in \Omega_{1}.
\end{equation*}
Let
\begin{equation*}
  u_2(x)=u_1(x)-Q_{k+1}(x',x_n-P_{\Omega}(x')),
\end{equation*}
where $P_{\Omega}\in \mathcal{P}_l$ corresponds to $\partial \Omega$ at $0$ (note that $\partial \Omega\in C^{l,\alpha}(0)$). Since $Q_{k+1}$ is a $(k+1)$-form and $P_{\Omega}(0)=|DP_{\Omega}(0)|=0$,
\begin{equation*}
D^{k+1}(Q_{k+1}(x',x_n-P_{\Omega}(x')))(0)=D^{k+1}(Q_{k+1}(x))(0).
\end{equation*}
Hence, $u_2(0)=\cdots= |D^{k+1}u_2(0)|=0$. In addition, $u_2$ satisfies
\begin{equation*}
\left\{\begin{aligned}
&\Delta u_2=f_2&& ~~\mbox{in}~~\Omega_1;\\
&u_2=g_2&& ~~\mbox{on}~~(\partial \Omega)_1,
\end{aligned}\right.
\end{equation*}
where $f_2\in C^{k+l-2,\alpha}(0)$ and (note that $Q_{k+1}$ is a $(k+1)$-form and $\partial \Omega\in C^{l,\alpha}(0)$)
\begin{equation*}
  |g_2(x)|\leq |g_1(x)|+|Q_{k+1}(x',x_n-P_{\Omega}(x'))|\leq C|x|^{k+l+\alpha},~~~\forall ~x\in (\partial \Omega)_1.
\end{equation*}

By \Cref{l-2.2} again, $u_2\in C^{k+2,\alpha}(0)$ and there exists a $(k+2)$-form $Q_{k+2}$ such that
\begin{equation*}
|u_2(x)-Q_{k+2}(x)|\leq C |x|^{k+2+\alpha},~~\forall ~x\in \Omega_{1}.
\end{equation*}
Let
\begin{equation}\label{e2.12}
  \begin{aligned}
u_3(x)&=u_2(x)-Q_{k+2}(x',x_n-P_{\Omega}(x'))\\
  &=u_1(x)-Q_{k+1}(x',x_n-P_{\Omega}(x'))-Q_{k+2}(x',x_n-P_{\Omega}(x')).
  \end{aligned}
\end{equation}
Then $u_3(0)=\cdots=|D^{k+2}u_3(0)|=0$. In addition, $u_3$ is a viscosity solution of
\begin{equation*}
\left\{\begin{aligned}
&\Delta u_3=f_3&& ~~\mbox{in}~~\Omega_1;\\
&u_3=g_3&& ~~\mbox{on}~~(\partial \Omega)_1,
\end{aligned}\right.
\end{equation*}
where $f_3\in C^{k+l-2,\alpha}(0)$ and
\begin{equation*}
  |g_3(x)|\leq |g_2(x)|+|Q_{k+2}(x',x_n-P_{\Omega}(x'))|\leq C|x|^{k+l+\alpha},~~~\forall ~x\in (\partial \Omega)_1,
\end{equation*}
where .

By \Cref{l-2.2} again, $u_3\in C^{k+2,\alpha}(0)$ and hence $u\in C^{k+2,\alpha}(0)$. By similar arguments again and again, $u\in C^{k+l,\alpha}(0)$ eventually and \cref{e1.2} holds. Therefore, the proof of \Cref{th1.1} is completed.\qed~\\

\begin{remark}\label{re2.2}
By checking the proof, we know that the polynomial $P$ can be written as \cref{e1.5} (see \cref{e2.12}).
\end{remark}
~\\

Next, we give the\\
\noindent\textbf{Proof of \Cref{co1.1}.} For $k=1$, \Cref{co1.1} reduces to \Cref{l-2.1}. For $k\geq 2$, $u\in C^{1,\alpha}(0)$ of course. Let
\begin{equation*}
\tilde{u}(x)=u(x)-u(0)-\sum_{i=1}^{n-1}P_{g,i}(0)x_i-u_n(0)(x_n-P_{\Omega}(x')).
\end{equation*}
Note that
\begin{equation*}
P_{\Omega}(0)=|DP_{\Omega}(0)|=0~~\mbox{ and }~~u_i(0)=P_{g,i}(0),~\forall ~1\leq i\leq n-1.
\end{equation*}
Hence, $\tilde{u}(0)=|D\tilde{u}(0)|=0$. In addition, $\tilde{u}$ satisfies
\begin{equation*}
\left\{\begin{aligned}
&\Delta \tilde{u}=\tilde{f}&& ~~\mbox{in}~~\Omega_1;\\
&\tilde{u}=\tilde{g}&& ~~\mbox{on}~~(\partial \Omega)_1,
\end{aligned}\right.
\end{equation*}
where $\tilde{f}\in C^{k-2,\alpha}(0)$, $\tilde{g}\in C^{k,\alpha}(0)$ and
\begin{equation*}
  \begin{aligned}
  |\tilde{g}(x)|&=|g(x)-g(0)-\sum_{i=1}^{n-1}P_{g,i}(0)x_i-u_n(0)(x_n-P_{\Omega}(x'))|\\
&=|g(x)-g(0)-\sum_{i=1}^{n}P_{g,i}(0)x_i+(P_{g,n}(0)-u_n(0))(x_n-P_{\Omega}(x'))
+P_{g,n}(0)P_{\Omega}(x')|\\
&\leq C|x|^{2},~~~\forall ~x\in (\partial \Omega)_1.
  \end{aligned}
\end{equation*}
Thus, $\tilde{g}(0)=|D\tilde{g}(0)|=0$. By \Cref{th1.1}, $\tilde{u}$ and hence $u\in C^{k,\alpha}(0)$.\qed~\\

Finally, we prove the higher regularity of free boundaries with the aid of the boundary pointwise regularity.\\
\noindent\textbf{Proof of \Cref{th1.2}.} Assume that
\begin{equation*}
  \partial \Omega\cap B_1=\left\{(x',x_n)\big| x_n=\varphi(x')\right\},
\end{equation*}
where $\varphi \in C^{1,\alpha}(T_1)$ and $\varphi(0)=|D\varphi(0)|=0$. Since $u=|Du|=0$ on $(\partial \Omega)_1$ and $(\partial \Omega)_1\in C^{1,\alpha}$, by \Cref{th1.1}, $u\in C^{2,\alpha}(x_0)$ for any $x_0\in (\partial \Omega)_1$. By combining with the interior regularity, $u\in C^{2,\alpha}(\bar\Omega')$ for any $\Omega'\subset\subset \bar\Omega\cap B_1$.

From $|D\varphi(0)|=0$ and $u=|Du|=0$ on $(\partial \Omega)_1$ again, $u_{ij}(0)=0$ for $i+j<2n$. Hence, $u_{nn}(0)=1$ by the equation $\Delta u=1$. Let $v(x)=u(x)-x_n^2/2$ and then $v$ satisfies
\begin{equation*}
\left\{\begin{aligned}
&   \Delta v=0~~~~\mbox{in}~~\Omega \cap B_1;\\
&   v=g~~~~\mbox{on}~~\partial \Omega\cap B_1;\\
&   v(0)=|Dv(0)|=|D^2v(0)|=0,
\end{aligned}\right.
\end{equation*}
where
\begin{equation*}
|g(x)|=|\frac{1}{2}x_n^2|\leq \|\partial \Omega\cap B_1\|_{C^{1,\alpha}(0)}|x|^{2+2\alpha},~\forall ~x\in (\partial \Omega)_1.
\end{equation*}
That is, $g\in C^{2,2\alpha}(0)$ and $g(0)=|Dg(0)|=|D^2g(0)|=0$. By \Cref{th1.1}, $v\in C^{2,2\alpha}(0)$ and hence $u\in C^{2,2\alpha}(0)$. Similarly, for any $x_0\in \partial \Omega\cap B_1$, $u\in C^{2,2\alpha}(x_0)$. Hence, $u\in C^{2,2\alpha}(\bar\Omega')$ for any $\Omega'\subset\subset \bar\Omega\cap B_1$.

Since $u_{nn}(0)=1$, $u_{nn}\geq 1/2$ in $\Omega_r$ for some $r>0$. Then $u_i/u_n\in C^{2\alpha} (\bar\Omega\cap B_{r})$. Note that $\varphi_i=u_i/u_n$. Thus, $\varphi\in C^{1,2\alpha}(T_r)$, i.e., $(\partial \Omega)_r\in C^{1,2\alpha}$. By considering other $x_0\in (\partial \Omega)_1$ similarly, we have $(\partial\Omega)_1\in C^{1,2\alpha}$.

Consider $v$ again and $g\in C^{2,4\alpha}(0)$ now (since $(\partial\Omega)_1\in C^{1,2\alpha}$). From \Cref{th1.1}, $v\in C^{2,4\alpha}(0)$. By similar arguments as above, $u\in C^{2,4\alpha}(\bar\Omega')$ for any $\Omega'\subset\subset \bar{\Omega}\cap B_1$. Therefore, $(\partial \Omega)_1\in C^{1, 4\alpha}$. Consider $v$ again and again and we have $(\partial \Omega)_1\in C^{2,\tilde{\alpha}}$ for some $0<\tilde{\alpha}<1$ eventually.

Let $v(x)=u(x)-(x_n-P_{\Omega}(x'))^2/2$ where $P_{\Omega}\in \mathcal{P}_2$ is the polynomial corresponding to $\partial \Omega$ at $0$ since $(\partial \Omega)_1\in C^{2,\tilde{\alpha}}$. Then $v$ satisfies
\begin{equation*}
\left\{\begin{aligned}
&   \Delta v=f~~~~\mbox{in}~~\Omega \cap B_1;\\
&   v=g~~~~\mbox{on}~~\partial \Omega\cap B_1;\\
&   v(0)=|Dv(0)|=|D^2v(0)|=0,
\end{aligned}\right.
\end{equation*}
where $f\in \mathcal{P}_2$ and
\begin{equation*}
|g(x)|=|(x_n-P_{\Omega}(x'))^2/2|\leq C|x|^{4+2\tilde{\alpha}},~\forall ~x\in (\partial \Omega)_1.
\end{equation*}
As before, by \Cref{th1.1}, $v\in C^{4, \tilde{\alpha}}(0)$ and hence $u\in C^{4, \tilde{\alpha}}(0)$. Similar to previous arguments, $u\in C^{4,\tilde{\alpha}}(\bar\Omega')$ for any $\Omega'\subset\subset \bar{\Omega}\cap B_1$ and then $(\partial \Omega)_1\in C^{3,\tilde{\alpha}}$.

Let $v(x)=u(x)-(x_n-P_{\Omega}(x'))^2/2$ where $P_{\Omega}\in \mathcal{P}_3$ now. Repeat above arguments and we have $u\in C^{\infty}(\bar\Omega\cap B_1)$ and $(\partial \Omega)_1\in C^{\infty}$ eventually. \qed~\\

\begin{remark}\label{re2.3}
In fact, we need a variation of \Cref{th1.1} in the proof. That is, if $g\in C^{k+l,\alpha}(0)$ is replaced by $g\in C^{k+\tilde{l},\tilde\alpha}(0)$ with $\tilde{l}+\tilde{\alpha}\leq l+\alpha$, we have $u\in C^{k+\tilde{l},\tilde\alpha}(0)$. This variation can be proved by almost the same proof. For clarity, we only give \Cref{th1.1} with $\tilde{l}=l$ and $\tilde{\alpha}=\alpha$.
\end{remark}

\begin{remark}\label{re2.4}
Maybe a more natural idea of proving \Cref{th1.2} is to consider $Du$ instead of $u-x_n^2/2$. The $u\in C^{2,\alpha}$ is easy to obtain. If $Du\in C^{2,\alpha}$ as well, $u\in C^{3,\alpha}$ and then $(\partial \Omega)_1\in C^{2,\alpha}$. By a series of iteration, the proof is completed.

For $1\leq i\leq n-1$, $u_i=0$ on $(\partial \Omega)_1$ and $Du_i(0)=0$. Hence, $u_i\in C^{2,\alpha}(0)$ by \Cref{th1.1}. However, we can't obtain $u_n\in C^{2,\alpha}(0)$ since $Du_n(0)=e_n\neq 0$. In addition, we can't obtain $u_i\in C^{2,\alpha}(x_0)$ for any $x_0\in (\partial \Omega)_1$ since $Du_i(x_0)\neq 0$ for other $x_0$.
\end{remark}

\section{Oblique derivative problem and application to the one-phase problem}\label{S3}
In this section, we give the detailed proofs of \Cref{th1.3-1} to \Cref{th1.4}. As in the proofs for Dirichlet problems, we intend to use compactness method to prove the regularity of solutions. Hence, we need to build a uniform estimate for solutions. First, we prove a Harnack type inequality.

\begin{lemma}\label{le3.1}
Let $u\geq 0$ be a viscosity solution of
\begin{equation*}
\left\{\begin{aligned}
&\Delta u=f&& ~~\mbox{in}~~\Omega_1;\\
&\beta\cdot Du=g&& ~~\mbox{on}~~(\partial \Omega)_1.
\end{aligned}\right.
\end{equation*}
Suppose that $\underset{B_1}{\mathrm{osc}}~\partial\Omega \leq \delta\leq \rho/8$, where $0<\rho<1$ depends only on $n$ and $a_0$.

Then for any $4\delta/\rho \leq R\leq 1/2$, we have
\begin{equation}\label{e3.2}
\sup_{\tilde{G}(R)}u\leq C\inf_{G(R/2)}u+CR\left(\|f\|_{L^{n}(\Omega_1)}
+\|g\|_{L^{\infty}((\partial \Omega)_1)}\right),
\end{equation}
where $C$ depends only on $n$ and $a_0$,
\begin{equation*}
G(R):=\left\{x\in\Omega \big||x'|<R,-\rho R<x_n<\rho R\right\}
\end{equation*}
and
\begin{equation*}
\tilde{G}(R):=\{x\in\Omega\big||x'|<R,x_n=\rho R\}.
\end{equation*}
\end{lemma}
\proof By the interior Harnack inequality,
\begin{equation*}
\inf_{\tilde{G}(R)}u\geq C\sup_{\tilde{G}(R)}u-R\|f\|_{L^{n}(\Omega_1)},
\end{equation*}
where $C$ depends only on $n$ and $\rho$. Hence, to prove \cref{e3.2}, we only need to show
\begin{equation}\label{e3.2-2}
\inf_{\tilde{G}(R)}u\leq C\inf_{G(R/2)}u+CR\left(\|f\|_{L^{n}(\Omega_1)}
+\|g\|_{L^{\infty}((\partial \Omega)_1)}\right).
\end{equation}
Without loss of generality, we assume that $\inf_{\tilde{G}(R)}u=1$.

Let
\begin{equation}\label{e3.14}
  v(x)=\frac{1}{2}+\frac{1}{4}\left(\left(\frac{x_n}{\rho R}\right)^2+\frac{x_n}{\rho R}-
  \frac{4|x'|^2}{R^2}\right).
\end{equation}
Then it can be verified easily that (by taking $\rho$ small enough)
\begin{equation*}
\left\{\begin{aligned}
&\Delta v\geq 0&& ~~\mbox{in}~~G(R);\\
&v\leq 1&&~~\mbox{on}~~ \tilde{G}(R);\\
&v\leq 0&&~~\mbox{on}~~\partial G(R)\backslash \left(\tilde{G}(R)\cup(\partial \Omega)_1\right);\\
&\beta\cdot Dv\geq 0&& ~~\mbox{on}~~(\partial \Omega)_1\cap G(R).
\end{aligned}\right.
\end{equation*}
Indeed, only the last inequality require some calculation. Since $\beta_n\geq a_0$, $\|\beta\|_{L^{\infty}}\leq 1$, $R\geq 4\delta/\rho$ and $\underset{B_1}{\mathrm{osc}}~\partial\Omega \leq \delta$, we have
\begin{equation*}
  \begin{aligned}
    \beta\cdot Dv&=\frac{\beta_n}{4\rho R}\left(\frac{2x_n}{\rho R}+1\right)
    -\sum_{i=1}^{n-1}\frac{2\beta_ix_i}{R^2}\\
    &\geq \frac{a_0}{4\rho R}\left(-\frac{1}{2}+1\right)-\frac{2(n-1)}{R}\\
    &\geq 0~~\mbox{ on }~~(\partial \Omega)_1\cap G(R),
  \end{aligned}
\end{equation*}
provided $\rho\leq a_0/(16n)$.

Let $w=u-v$ and then
\begin{equation*}
\left\{\begin{aligned}
&\Delta w\leq f&& ~~\mbox{in}~~G(R);\\
&w\geq 0&&~~\mbox{on}~~\partial G(R)\backslash (\partial \Omega)_1;\\
&\beta\cdot Dw\leq g&& ~~\mbox{on}~~(\partial \Omega)_1\cap G(R).
\end{aligned}\right.
\end{equation*}

By the Alexandrov-Bakel’man-Pucci maximum principle for oblique derivative problems (see \cite[Theorem 2.1]{MR3780142}),
\begin{equation*}
  w\geq -CR\|g\|_{L^{\infty}((\partial \Omega)_1)}-CR\|f\|_{L^{n}(\Omega_1)}
  ~~\mbox{ in }~~G(R),
\end{equation*}
where $C$ depends only on $n$ and $a_0$. Hence, by noting
\begin{equation*}
v\geq 1/8~~\mbox{ in }~~G(R/2),
\end{equation*}
we have
\begin{equation*}
  u=v+w\geq 1/8-CR\left(\|f\|_{L^{n}(\Omega_1)}
  +\|g\|_{L^{\infty}((\partial \Omega)_1)}\right)~~\mbox{ in }~~G(R/2).
\end{equation*}
That is, \cref{e3.2-2} holds.~\qed~\\

\begin{remark}\label{re3.1}
For the Dirichlet problem, we can obtain the equicontinuity up to the boundary by constructing proper barriers. In contrast, for the oblique derivative problem, we use the Harnack type inequality to show the equicontinuity since the solutions satisfies an equation on the boundary, which indicates that we should adopt the technique for the interior regularity rather than the technique for the boundary regularity (e.g. constructing barrier functions).
\end{remark}

\begin{remark}\label{re3.2}
Note that \cref{e3.2} is not a true Harnack inequality since it holds for $4\delta/\rho \leq R\leq 1/2$ other than $0<R\leq 1/2$. However, it is sufficient to provide the compactness in the proof (see \Cref{le3.3} below) and requires less smoothness of $\partial \Omega$.
\end{remark}

\begin{remark}\label{re3.3}
The construction of the auxiliary function $v$ is motivated by \cite{MR906819} (see Lemma 2.1 there) and has been used in \cite{MR3780142} (see Theorem 2.2 there).
\end{remark}
~\\

By a standard iteration argument, \Cref{le3.1} implies the following uniform estimate.
\begin{corollary}\label{co3.1}
Let $u$ be a viscosity solution of
\begin{equation*}
\left\{\begin{aligned}
&\Delta u=f&& ~~\mbox{in}~~\Omega_1;\\
&\beta\cdot Du=g&& ~~\mbox{on}~~(\partial \Omega)_1.
\end{aligned}\right.
\end{equation*}
Suppose that $\|u\|_{L^{\infty}(\Omega_1)}\leq 1$, $\|f\|_{L^{n}(\Omega_1)}\leq 1$, $\|g\|_{L^{\infty}((\partial \Omega)_1)}\leq 1$ and $\underset{B_1}{\mathrm{osc}}~\partial\Omega \leq \delta\leq \rho/8$, where $\rho$ is as in \Cref{le3.1}.

Then for $4\delta/\rho\leq R\leq 1/2$,
\begin{equation}\label{e3.3}
\underset{G(R)}{\mathrm{osc}}u\leq CR^{\alpha},
\end{equation}
where $0<\alpha<1$ is a universal constant and $C$ depends only on $n$ and $a_0$.
\end{corollary} ~~\\

Based on the above uniform estimate, we have the following equicontinuity for solutions.
\begin{lemma}\label{le3.2}
For any $\Omega'\subset\subset \bar{\Omega}\cap B_1$ and $\varepsilon>0$, there exists $\delta>0$ (depending only on $n,a_0,\varepsilon$ and $\Omega'$) such that if $u$ is a viscosity solution of
\begin{equation*}
\left\{\begin{aligned}
&\Delta u=f&& ~~\mbox{in}~~\Omega_1;\\
&\beta\cdot u=g&& ~~\mbox{on}~~(\partial \Omega)_1
\end{aligned}\right.
\end{equation*}
with $\|u\|_{L^{\infty}(\Omega_1)}\leq 1$, $\|f\|_{L^{n}(\Omega_1)}\leq 1$, $\|g\|_{L^{\infty}((\partial \Omega)_1)}\leq 1$ and $\underset{B_1}{\mathrm{osc}}~\partial\Omega \leq \delta$, then for any $x,y\in \Omega'$ with $|x-y|\leq \delta$,
\begin{equation*}
  |u(x)-u(y)|\leq \varepsilon.
\end{equation*}
\end{lemma}

\proof For any $\Omega'\subset\subset \bar{\Omega}\cap B_1$, $\varepsilon>0$ and $x,y\in \Omega'$, let $\delta>0$ to be specified later. Take $x_0\in (\partial \Omega)_1$ such that $|x-x_0|=\mathrm{dist}(x,(\partial \Omega)_1)$. By \Cref{co3.1} (a scaling version in fact), there exists $\delta_1>0$ (small enough) depending only on $n,a_0,\varepsilon$ and $\Omega'$ such that if $\underset{B_1}{\mathrm{osc}}~\partial\Omega \leq \rho\delta_1/4$, $|x-x_0|\leq \delta_1$ and $|y-x_0|\leq \delta_1$, we have
\begin{equation}\label{e.l3-2.1}
|u(x)-u(y)|\leq 2\underset{B(x_0,\delta_1)}{\mathrm{osc}}u\leq C\delta_1^{\alpha}\leq \varepsilon/2.
\end{equation}

If $|x-x_0|> \delta_1$, by the interior Lipschitz estimate for harmonic functions,
\begin{equation}\label{e.l3-2.2}
|u(x)-u(y)|\leq C\frac{|x-y|}{\delta_1},
\end{equation}
where $C$ depends only on $n$.

Take $\delta$ small enough such that $\delta\leq \rho\delta_1/4$ and $C\delta/\delta_1\leq \varepsilon/2$. Then by combining \cref{e.l3-2.1} and \cref{e.l3-2.2}, we derive the conclusion.~\qed~\\

In the following, we prove the boundary pointwise regularity for oblique derivative problems. First, we prove a key step.
\begin{lemma}\label{le3.3}
Let $0<\alpha<1$ and $u$ be a viscosity solution of
\begin{equation*}
\left\{\begin{aligned}
&\Delta u=f&& ~~\mbox{in}~~\Omega_1;\\
&\beta\cdot Du=g&& ~~\mbox{on}~~(\partial \Omega)_1.
\end{aligned}\right.
\end{equation*}
Suppose that $\|u\|_{L^{\infty}(\Omega_1)}\leq 1$, $\|f\|_{L^{n}(\Omega_1)}\leq\delta$, $\|g\|_{L^{\infty}((\partial \Omega)_1)}\leq \delta$, $\|\beta-\beta(0)\|_{L^{\infty}((\partial \Omega)_1)}\leq \delta$ and $\underset{B_1}{\mathrm{osc}}~\partial\Omega \leq \delta$, where $0<\delta<1$ depends only on $n,a_0$ and $\alpha$.

Then there exists $P\in \mathcal{P}_1$ such that
\begin{equation*}
\begin{aligned}
&\|u-P\|_{L^{\infty}(\Omega_{\eta})}\leq \eta^{1+\alpha},\\
&\|P\|\leq C_0,\\
&\beta(0)\cdot DP=0,\\
\end{aligned}
\end{equation*}
where $C_0$ depends only on $n$ and $a_0$, and $\eta$ depends also on $\alpha$.
\end{lemma}
\proof We prove the lemma by contradiction. Suppose that the lemma is false. Then there exist $0<\alpha<1$ and sequences of $u_m,f_m,g_m,\beta_m,\Omega_m$ such that
\begin{equation*}
\left\{\begin{aligned}
&\Delta u_m=f_m&& ~~\mbox{in}~~\Omega_m\cap B_1;\\
&\beta_m\cdot Du_m=g_m&& ~~\mbox{on}~~\partial \Omega_m\cap B_1
\end{aligned}\right.
\end{equation*}
with $\|u_m\|_{L^{\infty}(\Omega_m\cap B_1)}\leq 1$,
$\|f_m\|_{L^{n}(\Omega_m\cap B_1)}\leq 1/m$, $\|g_m\|_{L^{\infty}(\partial \Omega_m\cap B_1)}\leq 1/m$, $\|\beta_m-\beta_m(0)\|_{L^{\infty}((\partial \Omega)_1)}\leq 1/m$
and $\underset{B_1}{\mathrm{osc}}~\partial\Omega_m\leq 1/m$. In addition, for any $P\in \mathcal{P}_1$ with $\|P\|\leq C_0$ and $\beta_m(0)\cdot DP=0$,
\begin{equation}\label{el3.1}
  \|u_m-P\|_{L^{\infty}(\Omega_{m}\cap B_{\eta})}> \eta^{1+\alpha},
\end{equation}
where $C_0$ is to be specified later and $0<\eta<1$ is taken small such that
\begin{equation}\label{el3.2}
C_0\eta^{1-\alpha}<1/2.
\end{equation}

Note that $u_m$ are uniformly bounded ($\|u_m\|_{L^{\infty}(\Omega_m\cap B_1)}\leq 1$). Moreover, by \Cref{le3.2}, $u_m$ are equicontinuous. Precisely, for any $\Omega'\subset\subset B_1^+\cup T_1$, $\varepsilon>0$, there exist $\delta>0$ and $m_0$ such that for any $m\geq m_0$ and $x,y\in \Omega'\cap \bar{\Omega}_m$ with $|x-y|<\delta$, $|u(x)-u(y)|\leq \varepsilon$. Hence, there exists a subsequence (denoted by $u_m$ again) such that $u_m$ converges uniformly to some continuous function $\tilde{u}$ on compact subsets of $B_1^+\cup T_1$. Furthermore, there exists $\beta^0$ with $\beta^0_n\geq a_0$ such that $\beta_m(0)\rightarrow \beta^0$. By the closedness of viscosity solutions (e.g., see \cite[Proposition 2.1.]{MR3780142}), $\tilde{u}$ is a viscosity solution of
\begin{equation*}
\left\{\begin{aligned}
 &\Delta \tilde{u}=0&& ~~\mbox{in}~~B_{1}^+;\\
 &\beta^0\cdot D\tilde{u}=0&&~~\mbox{on}~~T_{1}.
\end{aligned}\right.
\end{equation*}

By the boundary regularity for homogeneous equations on flat boundaries
(e.g. see \cite[Theorem 4.1 and Theorem 4.2]{MR3780142}), there exists $\tilde P\in \mathcal{P}_1$ such that
\begin{equation}\label{el3.3}
\begin{aligned}
  &|\tilde{u}(x)-\tilde{P}(x)|\leq C_0 |x|^{2}, ~~\forall ~x\in B_{1/2}^+,\\
  &\|\tilde{P}\|\leq C_0/2,\\
  &\beta^{0}\cdot D\tilde{P}=0.
\end{aligned}
\end{equation}

Since $\beta_m(0)\rightarrow \beta^0$, there exists $P_m\in \mathcal{HP}_1$ such that $\beta_m(0)\cdot (D\tilde{P}+DP_m)=0$ and $\|P_m\|\rightarrow 0$ as $m\rightarrow \infty$. Thus, \cref{el3.1} holds for $\tilde{P}+P_m$. Let $m\rightarrow \infty$ in \cref{el3.1} and we have
\begin{equation*}
    \|\tilde{u}-\tilde{P}\|_{L^{\infty}(B_{\eta}^+)}\geq \eta^{1+\alpha}.
\end{equation*}
On the other hand, from \cref{el3.2} and \cref{el3.3},
\begin{equation*}
  \|\tilde{u}-\tilde{P}\|_{L^{\infty}(B_{\eta}^+)}\leq \eta^{1+\alpha}/2,
\end{equation*}
which is a contradiction. ~\qed~\\

\begin{remark}\label{re3.5}
Usually, to prove \Cref{le3.3}, we solve an equation to approximate $u$ (see the proofs of Lemma 5.1 and Lemma 6.3 in \cite{MR3780142}). Instead, the compactness method avoids the solvability. This is one of the main advantages of the method of compactness.
\end{remark}

\begin{remark}\label{re3.6}
Note that \Cref{le3.2} is not a true equicontinuity up to the boundary. However, it is enough to provide the compactness in the proof of \Cref{le3.3}. The benefit is that we don't require the smoothness of $\partial \Omega\cap B_1$ and hence we can develop the pointwise regularity.
\end{remark}
~\\

Now, we can prove the boundary pointwise $C^{1,\alpha}$ regularity.
\begin{theorem}\label{th3.1}
Let $0<\alpha <1$ and $u$ be a viscosity solution of
  \begin{equation*}
    \left\{
    \begin{aligned}
      &\Delta u=f~~ &&\mbox{in}~~ \Omega_1; \\
      &\beta  \cdot Du = g~~ &&\mbox{on}~~(\partial \Omega)_1.
    \end{aligned}
    \right.
  \end{equation*}
Suppose that $f\in C^{-1,\alpha}(0)$, $g\in C^{\alpha}(0)$, $\beta\in C^{\alpha}(0)$ and $[\partial \Omega\cap B_1]_{C^{0,1}(0)} \leq \delta$, where $\delta$ is as in \Cref{le3.3}.

Then $u$ is $C^{1,\alpha}$ at $0$, i.e., there exists $P\in \mathcal{P}_1$ such that
\begin{equation}\label{et3.1}
\begin{aligned}
&|u(x)-P(x)|\leq C|x|^{1+\alpha}\left(\|u\|_{L^{\infty}(\Omega_1)}
+\|f\|_{C^{-1,\alpha}(0)}+\|g\|_{C^{\alpha}(0)}\right), ~~\forall x\in \Omega_1,\\
&|Du(0)|\leq C\left(\|u\|_{L^{\infty}(\Omega_1)}
+\|f\|_{C^{-1,\alpha}(0)}+\|g\|_{C^{\alpha}(0)}\right),
\end{aligned}
\end{equation}
and
\begin{equation}\label{et3.2}
\beta(0)\cdot DP=0,
\end{equation}
where $C$ depends only on $n,a_0,\alpha$ and $\|\beta \|_{C^{\alpha }(0)}$.
\end{theorem}
\proof We assume that $g(0)=0$. Otherwise, we may consider $\tilde{u}=u-g(0)x_n/\beta_n(0)$. Let $\delta$ be as in \Cref{le3.3}, which depends only on $n,a_0$ and $\alpha$. Without loss of generality, we assume that
\begin{equation}\label{e.tC1a-ass}
  \begin{aligned}
    &\|u\|_{L^{\infty}(\Omega_1)}\leq 1,~ \|f\|_{C^{-1,\alpha}(0)}\leq \delta,
    ~[g]_{C^{\alpha}(0)}\leq \frac{\delta}{2},
    ~[\beta]_{C^{\alpha}(0)}\leq \frac{\delta}{2C_1},\\
  \end{aligned}
\end{equation}
where $C_1$ is a constant (depending only on $n,a_0$ and $\alpha$) to be specified later.

To show $u\in C^{1,\alpha}(0)$, we only need to prove that there exists a sequence of $P_m\in \mathcal{P}_1$ ($m\geq -1$) such that for all $m\geq 0$,
\begin{equation}\label{e3.16}
\|u-P_m\|_{L^{\infty }(\Omega _{\eta^{m}})}\leq \eta ^{m(1+\alpha )},
\end{equation}
\begin{equation}\label{e3.17}
|P_m(0)-P_{m-1}(0)|+\eta^{m}|DP_m-DP_{m-1}|\leq C_0\eta ^{m(1+\alpha)}
\end{equation}
and
\begin{equation}\label{e3.18}
\beta(0)\cdot DP_m=0,
\end{equation}
where $C_0$ and $\eta$ are constants as in \Cref{le3.3}.

We prove the above by induction. For $m=0$, by setting $P_0\equiv P_{-1}\equiv 0$, the conclusion holds clearly. Suppose that the conclusion holds for $m$. We need to prove that the conclusion holds for $m+1$.

Let $r=\eta ^{m}$, $y=x/r$ and
\begin{equation}\label{e-v1}
  \tilde{u}(y)=\frac{u(x)-P_{m}(x)}{r^{1+\alpha}}.
\end{equation}
Then $\tilde{u}$ satisfies
\begin{equation*}
\left\{\begin{aligned}
&\Delta \tilde{u}=\tilde{f}&& ~~\mbox{in}~~\tilde{\Omega}\cap B_1;\\
&\tilde{\beta}\cdot D\tilde{u}=\tilde{g}&& ~~\mbox{on}~~\partial \tilde{\Omega}\cap B_1,
\end{aligned}\right.
\end{equation*}
where
\begin{equation*}
  \tilde{f}(y)=\frac{f(x)}{r^{\alpha-1}},~\tilde{g}(y)=\frac{g(x)-\beta(x)\cdot DP_m}{r^{\alpha}},~
  \tilde{\beta}(y)=\beta(x)
  ~~\mbox{and}~~  \tilde{\Omega}=\frac{\Omega}{r}.
\end{equation*}

By \cref{e3.17}, there exists a constant $C_1$ depending only on $n,a_0$ and $\alpha$ such that $|DP_i|\leq C_1$ ($\forall~0\leq i\leq m$). Then it is easy to verify that
\begin{equation*}
\begin{aligned}
&\|\tilde{u}\|_{L^{\infty}(\tilde{\Omega}\cap B_1)}\leq 1, ~(\mathrm{by~} \cref{e3.16} ~\mathrm{ and }~ \cref{e-v1})\\
&\|\tilde{f}\|_{L^{n}(\tilde{\Omega}\cap B_1)}=\frac{\|f\|_{L^{n}(\Omega\cap B_r)}}{r^{\alpha}}\leq \delta, ~(\mathrm{by }~\cref{e.tC1a-ass})\\
&\|\tilde{g}\|_{L^{\infty}(\partial \tilde{\Omega}\cap B_1)}
\leq \frac{1}{r^{\alpha}}\left([g]_{C^{\alpha}(0)}r^{\alpha}
  +C_1[\beta]_{C^{\alpha}(0)}r^{\alpha}\right)\leq \delta,  ~(\mathrm{by}~ \cref{e.tC1a-ass} ~\mathrm{ and }~ \cref{e3.18})\\
&\|\tilde{\beta}-\tilde{\beta}(0)\|_{L^{\infty}(\partial \tilde{\Omega}\cap B_1)}
=\|\beta-\beta(0)\|_{L^{\infty}(\partial \tilde{\Omega}\cap B_r)}\leq [\beta]_{C^{\alpha}(0)}r^{\alpha}  \leq \delta ,~(\mathrm{by}~ \cref{e.tC1a-ass})\\
&\underset{B_1}{\mathrm{osc}}~\partial\tilde{\Omega}=
\frac{1}{r}\underset{B_r}{\mathrm{osc}}~\partial\Omega \leq [\partial \Omega\cap B_1]_{C^{0,1}(0)}\leq \delta.
\end{aligned}
\end{equation*}

From \Cref{le3.3}, there exists $\tilde{P}\in \mathcal{P}_1$ such that
\begin{equation*}
\begin{aligned}
&\|\tilde{u}-\tilde{P}\|_{L^{\infty }(\tilde{\Omega} _{\eta})}\leq \eta ^{1+\alpha},\\
&\|\tilde{P}\|\leq C_0,\\
&\tilde\beta(0)\cdot D\tilde{P}=0.
\end{aligned}
\end{equation*}
Let $P_{m+1}(x)=P_{m}(x)+r^{1+\alpha}\tilde{P}(y)$. Then \cref{e3.17} and \cref{e3.18} hold for $m+1$. Recalling \cref{e-v1}, we have
\begin{equation*}
  \begin{aligned}
&\|u-P_{m+1}\|_{L^{\infty}(\Omega_{\eta^{m+1}})}\\
&= \|u-P_{m}-r^{1+\alpha}\tilde{P}(y)\|_{L^{\infty}(\Omega_{\eta r})}\\
&= \|r^{1+\alpha}\tilde{u}-r^{1+\alpha}\tilde{P}\|_{L^{\infty}(\tilde{\Omega}_{\eta})}\\
&\leq r^{1+\alpha}\eta^{1+\alpha}\\
&=\eta^{(m+1)(1+\alpha)}.
  \end{aligned}
\end{equation*}
Hence, \cref{e3.16} holds for $m+1$. By induction, the proof is completed.~\qed~\\

Next, we prove a lemma similar to \Cref{l-2.2}.
\begin{lemma}\label{le-3.1}
Let $0<\alpha<1$ and $u$ be a viscosity solution of
\begin{equation*}
\left\{\begin{aligned}
&\Delta u=f&& ~~\mbox{in}~~\Omega_1;\\
&\beta\cdot Du=g&& ~~\mbox{on}~~(\partial \Omega)_1.
\end{aligned}\right.
\end{equation*}
Suppose that $u\in C^{k,\alpha}(0)(k\geq 1)$, $f\in C^{k-1,\alpha}(0)$, $g\in C^{k,\alpha}(0)$, $\beta\in C^{\alpha}(0)$ and $(\partial \Omega)_1\in C^{1,\alpha}(0)$. Moreover, assume that
\begin{equation*}
u(0)=\cdots=|D^ku(0)|=|Dg(0)|\cdots=|D^{k-1}g(0)|=0.
\end{equation*}

Then $u\in C^{k+1,\alpha}(0)$. That is, there exists $P\in \mathcal{HP}_{k+1}$ such that
\begin{equation*}
  \begin{aligned}
&|u(x)-P(x)|\leq C |x|^{k+1+\alpha}\left(\|u\|_{L^{\infty}}+\|f\|_{C^{k-1,\alpha}(0)}
+\|g\|_{C^{k,\alpha}(0)}\right),~~\forall ~x\in \Omega_{1},\\
&|D^{k+1}u(0)|\leq C\left(\|u\|_{L^{\infty}}+\|f\|_{C^{k-1,\alpha}(0)}
+\|g\|_{C^{k,\alpha}(0)}\right),\\
  \end{aligned}
\end{equation*}
and
\begin{equation*}
\Delta P\equiv P_f,~~\beta(0)\cdot DP(x',0)\equiv P_g(x',0),
\end{equation*}
where $C$ depends only on $n,k,a_0,\alpha$, $\|\beta\|_{C^{\alpha}(0)}$ and $\|(\partial \Omega)_1\|_{C^{1,\alpha}(0)}$.
\end{lemma}

As in \Cref{S2}, we prove above lemma by induction. For $k=1$, the lemma reduces to \Cref{th3.1}. Suppose that the lemma holds for $k\leq k_0-1$ and we need to prove the lemma for $k=k_0$. First, we prove the following lemma which is a key step towards the conclusion of \Cref{le-3.1}.

\begin{lemma}\label{le3.4}
Let $1\leq k\leq k_0$, $0<\alpha<1$ and $u\in C^{k,\alpha}(0)$ be a viscosity solution of
\begin{equation*}
\left\{\begin{aligned}
&\Delta u+P=f&& ~~\mbox{in}~~\Omega_1;\\
&\beta\cdot Du=g&& ~~\mbox{on}~~(\partial \Omega)_1,
\end{aligned}\right.
\end{equation*}
where $P\in \mathcal{HP}_{k-1}$. Suppose that
\begin{equation*}
  \begin{aligned}
    &\|u\|_{L^{\infty}(\Omega_1)}\leq 1, u(0)=\cdots=|D^ku(0)|=0,\\
    &|f(x)|\leq \delta|x|^{k-2+\alpha}, ~~\forall x\in \Omega_1,\\
    &|g(x)|\leq \delta|x|^{k-1+\alpha},~~\forall x\in (\partial \Omega)_1,\\
    &|\beta(x)-\beta(0)|\leq \delta|x|^{\alpha},~~\forall x\in (\partial \Omega)_1,\\
    &\|(\partial \Omega)_1\|_{C^{1,\alpha}(0)} \leq \delta,\\
    &\|P\|\leq 1,
  \end{aligned}
\end{equation*}
where $\delta>0$ depending only on $n,k,a_0$ and $\alpha$.

Then there exists $Q\in \mathcal{HP}_{k+1}$ such that
\begin{equation*}
  \begin{aligned}
    & \|u-Q\|_{L^{\infty}(\Omega_{\eta})}\leq \eta^{k+1+\alpha},\\
    &\|Q\|\leq C_0,\\
    &\Delta Q+P\equiv 0,\\
    &\beta(0)\cdot DQ(x',0)\equiv 0,\\
  \end{aligned}
\end{equation*}
where $C_0$ depends only on $n,a_0$ and $k$, and $\eta$ depends also on $\alpha$.
\end{lemma}
\begin{remark}\label{re3.4}
Note that $\beta(0)\cdot DQ(x',0)\equiv 0$ is equivalent to that $\beta(0)\cdot DQ(x)$ is a $k$-form.
\end{remark}
\proof As before, we prove the lemma by contradiction. Suppose that the conclusion is false. Then there exist $0<\alpha<1$ and sequences of $u_m,f_m,g_m,\beta_m,\Omega_m, P_m$ ($m\geq 1$) satisfying $u_m\in C^{k,\alpha}(0)$ and
\begin{equation*}
\left\{\begin{aligned}
&\Delta u_m+P_m=f_m&& ~~\mbox{in}~~\Omega_m\cap B_1;\\
&\beta_m\cdot Du_m=g_m&& ~~\mbox{on}~~\partial \Omega_m\cap B_1.
\end{aligned}\right.
\end{equation*}
In addition,
\begin{equation*}
  \begin{aligned}
    &\|u_m\|_{L^{\infty}(\Omega_m\cap B_1)}\leq 1, u_m(0)=\cdots=|D^ku_m(0)|=0,\\
    &|f_m(x)|\leq |x|^{k-2+\alpha}/m, ~~\forall x\in \Omega_1,\\
    &|g_m(x)|\leq |x|^{k-1+\alpha}/m,~~\forall x\in (\partial \Omega)_1,\\
    &|\beta_m-\beta_m(0)|\leq |x|^{\alpha}/m,~~\forall x\in (\partial \Omega)_1,\\
    &\|\partial \Omega_m\cap B_1\|_{C^{1,\alpha}(0)} \leq 1/m,\\
    &\|P_m\|\leq 1.
  \end{aligned}
\end{equation*}
But for any $Q\in \mathcal{HP}_{k+1}$ satisfying
\begin{equation*}
  \begin{aligned}
    &\|Q\|\leq C_0,\\
    &\Delta Q+P_m\equiv 0,\\
    &\beta_m(0)\cdot DQ(x',0)\equiv 0,\\
  \end{aligned}
\end{equation*}
we have
\begin{equation}\label{e3.4}
  \|u_m-Q\|_{L^{\infty}(\Omega_{m}\cap B_{\eta})}> \eta^{k+1+\alpha},
\end{equation}
where $C_0$ is to be specified later and $0<\eta<1$ is taken small such that
\begin{equation}\label{e3.5}
C_0\eta^{1-\alpha}<1/2.
\end{equation}

As in the proof of \Cref{le3.3}, $u_m$ are uniformly bounded and equicontinuous. Hence, there exist $\tilde u:B_1^+\cup T_1\rightarrow \mathbb{R}$, $\beta^0\in \mathbb{R}^n$ with $\beta^0_n\geq a_0$ and $\tilde P\in \mathcal{HP}_{k-1}$ such that $u_m\rightarrow\tilde u$ uniformly in compact subsets of $ B_1^+\cup T_1$, $\beta_m(0)\rightarrow \beta^0$, $P_m\rightarrow \tilde P$ and
\begin{equation*}
\left\{\begin{aligned}
&\Delta \tilde u+\tilde P=0&& ~~\mbox{in}~~B_{1}^+;\\
&\beta^0\cdot D\tilde u=0&& ~~\mbox{on}~~T_{1}.
\end{aligned}\right.
\end{equation*}

By the boundary $C^{k,\alpha}$ estimate for $u_m$ (\Cref{le-3.1} for $k-1$ since $k\leq k_0$) and noting $u_m(0)=\cdots=|D^{k}u_m(0)|=0$, we have
\begin{equation*}
\|u_m\|_{L^{\infty }(\Omega_m\cap B_r)}\leq Cr^{k+\alpha} ,~~~~\forall~0<r<1.
\end{equation*}
Since $u_m$ converges to $u$ uniformly,
\begin{equation*}
\|\tilde u\|_{L^{\infty }(B_r^+)}\leq Cr^{k+\alpha}, ~~~~\forall~0<r<1.
\end{equation*}
Hence, $\tilde u(0)=\cdots=|D^{k}\tilde u(0)|=0$. By the boundary estimate for $\tilde u$, there exists $\tilde{Q}\in \mathcal{HP}_{k+1}$ such that
\begin{equation}\label{e3.7}
\begin{aligned}
&|\tilde u(x)-\tilde{Q}(x)|\leq C_0 |x|^{k+2}, ~~\forall ~x\in B_{1}^+,\\
&\|\tilde{Q}\|\leq C_0/2,\\
&\Delta \tilde{Q}+\tilde P\equiv 0,\\
&\beta^0\cdot D\tilde{Q}(x',0)\equiv 0,\\
\end{aligned}
\end{equation}
where $C_0$ depends only on $n,a_0$ and $k$.

Since $\beta_m(0)\rightarrow \beta^0$ and $P_m\rightarrow \tilde{P}$, there exist $Q_m\in \mathcal{HP}_{k+1}$ such that $\|Q_m\|\rightarrow 0$ and
\begin{equation*}
\begin{aligned}
&\|\tilde{Q}+Q_m\|\leq C_0,\\
&\Delta (\tilde{Q}+ Q_m)+P_m\equiv 0,\\
&\beta_m(0)\cdot D(\tilde{Q}+Q_m)(x',)\equiv 0.
\end{aligned}
\end{equation*}
Thus, \cref{e3.4} holds for $Q=\tilde{Q}+Q_m$. Let $m\rightarrow \infty$ in \cref{e3.4} and we have
\begin{equation*}
    \|\tilde u-\tilde{Q}\|_{L^{\infty}(B_{\eta}^+)}\geq \eta^{k+1+\alpha},
\end{equation*}
However, by \cref{e3.5} and \cref{e3.7},
\begin{equation*}
  \|\tilde u-\tilde{Q}\|_{L^{\infty}(B_{\eta}^+)}\leq \eta^{k+1+\alpha}/2,
\end{equation*}
which is a contradiction.  ~\qed~\\

Now, we give the\\
\noindent\textbf{Proof of \Cref{le-3.1}.} Since we have assumed that \Cref{le-3.1} holds for $k_0-1$, $u\in C^{k_0,\alpha}(0)$. By induction, we only need to prove \Cref{le-3.1} for $k_0$, i.e., $u\in C^{k_0+1,\alpha}(0)$. Without loss of generality, by a proper transformation, we can assume as in the proof of \Cref{l-2.2} that
\begin{equation*}
\left\{\begin{aligned}
&\Delta u+P=f&& ~~\mbox{in}~~\Omega_1;\\
&u=g&& ~~\mbox{on}~~(\partial \Omega)_1
\end{aligned}\right.
\end{equation*}
for some $P\in\mathcal{HP}_{k_0-1}$ and
\begin{equation}\label{e3.8}
  \begin{aligned}
    &\|u\|_{L^{\infty}(\Omega_1)}\leq 1,u(0)=\cdots=|D^{k_0}u(0)|=0,\\
    &|f(x)|\leq \delta |x|^{k_0-1+\alpha}, ~~\forall x\in \Omega_1,\\
    &|g(x)|\leq \delta |x|^{k_0+\alpha}/3,~~\forall x\in (\partial \Omega)_1,\\
    &|\beta(x)-\beta(0)|\leq \delta|x|^{\alpha}/(3C_1),~~\forall x\in (\partial \Omega)_1,\\
    &\|(\partial \Omega)_1\|_{C^{1,\alpha}(0)} \leq \delta/(3C_1),\\
    &\|P\|\leq 1,
  \end{aligned}
\end{equation}
where $\delta$ is as in \Cref{le3.4} (with $k=k_0$) and $C_1$ depending only on $n,k_0,a_0$ and $\alpha$ is to be chosen later.

To prove \Cref{le-3.1} for $k_0$, we only to show that there exists a sequence of $Q_m\in \mathcal{HP}_{k_0+1}$ ($m\geq 0$) such that for all $m\geq 1$,
\begin{equation}\label{e3.9}
\|u-Q_m\|_{L^{\infty }(\Omega _{\eta^{m}})}\leq \eta ^{m(k_0+1+\alpha )},
\end{equation}
\begin{equation}\label{e3.12}
\|Q_m-Q_{m-1}\|\leq C_0\eta ^{m\alpha},
\end{equation}
and
\begin{equation}\label{e3.10}
\Delta Q_m+P\equiv 0,~~\beta(0)\cdot DQ_m(x',0)\equiv 0,
\end{equation}
where $C_0$ and $\eta$ are the constants as in \Cref{le3.4}.

We prove the above by induction. For $m=1$, by \Cref{le3.4} and setting $Q_0\equiv 0$, the conclusion holds clearly. Suppose that the conclusion holds for $m$. We need to prove that the conclusion holds for $m+1$.

Let $r=\eta^{m}$, $y=x/r$ and
\begin{equation}\label{e3.13}
  \tilde{u}(y)=\frac{u(x)-Q_m(x)}{r^{k_0+1+\alpha}}.
\end{equation}
Then $\tilde{u}$ satisfies
\begin{equation*}
\left\{\begin{aligned}
&\Delta \tilde{u}=\tilde{f}&& ~~\mbox{in}~~\tilde{\Omega}\cap B_1;\\
&\tilde{\beta}\cdot D\tilde{u}=\tilde{g}&& ~~\mbox{on}~~\partial \tilde{\Omega}\cap B_1,
\end{aligned}\right.
\end{equation*}
where
\begin{equation*}
 \tilde{f}(y)=\frac{f(x)}{r^{k_0-1+\alpha}},~~
 \tilde{g}(y)=\frac{g(x)-\beta(x)\cdot DQ_m(x)}{r^{k_0+\alpha}},~~
 \tilde{\beta}(y)=\beta(x),~~
 \tilde{\Omega}=\frac{\Omega}{r}.
\end{equation*}
By \cref{e3.12}, there exists $C_1$ depends only on $n,k_0,a_0$ and $\alpha$ such that $\|Q_i\|\leq C_1$ ($\forall~0\leq i\leq m$). Since $\beta(0)\cdot DQ_m$ is a $k_0$-form (see \cref{e3.10}) and $\partial \Omega\in C^{1,\alpha}(0)$,
\begin{equation}\label{e3.25}
\begin{aligned}
|\beta(0)\cdot DQ_m(x)|&\leq C_1|x|^{k_0-1}|x_n|
\leq C_1\|(\partial \Omega)_1\|_{C^{1,\alpha}(0)}|x|^{k_0+\alpha},~\forall ~x\in (\partial \Omega)_1.
\end{aligned}
\end{equation}
Then
\begin{equation*}
  \begin{aligned}
    & \|\tilde{u}\|_{L^{\infty}(\tilde{\Omega}\cap B_1)}\leq 1, ~(\mathrm{by}~ \cref{e3.9} ~\mbox{and}~ \cref{e3.13})\\
    & |\tilde{f}(y)|=\frac{|f(x)|}{r^{k_0-1+\alpha}}\leq \delta|y|^{k_0-1+\alpha},~~\forall y\in \tilde\Omega_1, ~(\mathrm{by}~ \cref{e3.8})\\
    &|\tilde{g}(y)|\leq \frac{1}{r^{k_0+\alpha}}\left(|g(x)|+|\beta(x)-\beta(0)||DQ_m(x)|+|\beta(0)\cdot DQ_m(x)|\right)\\
    &~~~~~~~\leq \frac{1}{r^{k_0+\alpha}}\left(\frac{\delta}{3}|x|^{k_0+\alpha}
      +\frac{\delta}{3C_1}\cdot C_1|x|^{k_0+\alpha}
      +C_1\cdot\frac{\delta}{3C_1}|x|^{k_0+\alpha}\right)\\
     &~~~~~~~\leq \delta |y|^{k_0+\alpha},~~\forall y\in (\partial \tilde\Omega)_1, ~(\mathrm{by}~ \cref{e3.8}~\mbox{and}~ \cref{e3.25})\\
     &|\tilde{\beta}(y)-\tilde{\beta}(0)|=|\beta(x)-\beta(0)|\leq r^{\alpha}[\beta]_{C^{\alpha}(0)}|y|^{\alpha}\leq \delta |y|^{\alpha},~~\forall y\in (\partial \tilde\Omega)_1, ~(\mathrm{by}~ \cref{e3.8})\\
     &\|\partial \tilde{\Omega}\cap B_1\|_{C^{1,\alpha}(0)}\leq r^{\alpha}\|(\partial \Omega)_1\|_{C^{1,\alpha}(0)} \leq \delta.   ~(\mathrm{by}~ \cref{e3.8})
  \end{aligned}
\end{equation*}

By \Cref{le3.4}, there exists $\tilde{Q}\in \mathcal{HP}_{k_0+1}$ such that
\begin{equation*}
\begin{aligned}
    &\|\tilde{u}-\tilde{Q}\|_{L^{\infty }(\tilde{\Omega} _{\eta})}\leq \eta ^{k_0+1+\alpha},\\
    &\|\tilde{Q}\|\leq C_0,\\
    &\Delta\tilde{Q}\equiv 0,\\
    &\beta(0)\cdot D\tilde{Q}(x',0)\equiv 0.
\end{aligned}
\end{equation*}
Let $Q_{m+1}(x)=Q_m(x)+r^{k_0+1+\alpha}\tilde{Q}(y)=Q_m(x)+r^{\alpha}\tilde{Q}(x)$.
Then \cref{e3.12} and \cref{e3.10} hold for $m+1$. By recalling \cref{e3.13}, we have
\begin{equation*}
  \begin{aligned}
&\|u-Q_{m+1}\|_{L^{\infty}(\Omega_{\eta^{m+1}})}\\
&= \|u-Q_m-r^{\alpha}\tilde{Q}\|_{L^{\infty}(\Omega_{\eta r})}\\
&= \|r^{k_0+1+\alpha}\tilde{u}-r^{k_0+1+\alpha}\tilde{Q}\|_{L^{\infty}(\tilde{\Omega}_{\eta})}\\
&\leq r^{k_0+1+\alpha}\eta^{k_0+1+\alpha}\\
&=\eta^{(m+1)(k_0+1+\alpha)}.
  \end{aligned}
\end{equation*}
Hence, \cref{e3.9} holds for $m+1$. By induction, the proof is completed.\qed~\\

Next, we prove \Cref{th1.3} with the aid of \Cref{le-3.1}.\\
\noindent\textbf{Proof of \Cref{th1.3}.} Throughout this proof, $C$ always denotes a constant depending only on $n,k,l,a_0,\alpha$, $\|\beta\|_{C^{l-1,\alpha}(0)}$ and $\|(\partial \Omega)_1\|_{C^{l,\alpha}(0)}$. Without loss of generality, we assume as before
\begin{equation*}
\|u\|_{L^{\infty}(\Omega_1)}
+\|f\|_{C^{k+l-2,\alpha}(0)}+\|g\|_{C^{k+l-1,\alpha}(0)}\leq 1.
\end{equation*}

Since $g\in C^{k+l-1,\alpha}(0)$ and $\beta\in C^{l-1,\alpha}(0)$,
\begin{equation}\label{e3.19}
  |g(x)-P_g(x)|\leq |x|^{k+l-1+\alpha},~~~\forall ~x\in (\partial \Omega)_1.
\end{equation}
and
\begin{equation}\label{e3.20}
  |\beta(x)-P_{\beta}(x)|\leq [\beta]_{C^{l-1,\alpha}(0)}|x|^{l-1+\alpha},~~~\forall ~x\in (\partial \Omega)_1.
\end{equation}
Note that $P_g\in \mathcal{HP}_{k}$ and $P_{\beta}$ is a vector valued polynomial. Take a polynomial $P_0\in \mathcal{HP}_{k+l}$ such that
\begin{equation*}
  P_{\beta}\cdot DP_0\equiv P_{g}.
\end{equation*}

Set $u_1=u-P_0$ and then $u_1\in C^{k,\alpha}(0)$ and $u_1(0)=\cdots=|D^ku_1(0)|=0$. In addition, $u_1$ is a viscosity solution of
\begin{equation*}
\left\{\begin{aligned}
&\Delta u_1=f_1&& ~~\mbox{in}~~\Omega_1;\\
&\beta\cdot Du_1=g_1&& ~~\mbox{on}~~(\partial \Omega)_1,
\end{aligned}\right.
\end{equation*}
where $f_1=f-\Delta P_0$ and
\begin{equation*}
g_1=g-\beta\cdot P_0=g-P_g-\left(\beta-P_{\beta}\right)\cdot DP_0.
\end{equation*}
Hence, by \cref{e3.19}, \cref{e3.20} and noting $P_0\in \mathcal{HP}_{k+l}$,
\begin{equation}\label{e3.11}
|g_1(x)|\leq C|x|^{k+l-1+\alpha},~~\forall ~x\in (\partial \Omega)_1.
\end{equation}
By \Cref{le-3.1}, $u_1\in C^{k+1,\alpha}(0)$. That is, there exists $P_{k+1}\in \mathcal{HP}_{k+1}$ such that $\beta(0)\cdot DP_{k+1}$ is a $k$-form and
\begin{equation*}
  |u_1(x)-P_{k+1}(x)|\leq C|x|^{k+1+\alpha},~~~\forall ~x\in \Omega_1.
\end{equation*}

Let
\begin{equation*}
  u_2(x)=u_1(x)-P_{k+1}(x',x_n-P_{\Omega}),
\end{equation*}
where $P_{\Omega}\in \mathcal{P}_{l}$ corresponds to $\partial \Omega$ at $0$ (note that $\partial \Omega\in C^{l,\alpha}(0)$). Note that $P_{\Omega}(0)=|DP_{\Omega}(0)|=0$ and then $u_2(0)=\cdots=|D^{k+1}u_2(0)|=0$. In addition, $u_2$ satisfies
\begin{equation*}
\left\{\begin{aligned}
&\Delta u_2=f_2&& ~~\mbox{in}~~\Omega_1;\\
&\beta\cdot Du_2=g_2&& ~~\mbox{on}~~(\partial \Omega)_1,
\end{aligned}\right.
\end{equation*}
where $f_2\in C^{k+l-2,\alpha}(0)$,
\begin{equation*}
  \begin{aligned}
g_2(x)=&g_1(x)-\beta(x)\cdot\left(DP_{k+1}(x',x_n-P_{\Omega}(x'))
-P_{k+1,n}(x',x_n-P_{\Omega}(x'))DP_{\Omega}(x')\right)\\
=&g_1(x)-\tilde{g}_2(x)-\bar{g}_2(x)-\hat{g}_2(x),\\
\end{aligned}
\end{equation*}
and
\begin{equation*}
\begin{aligned}
\tilde{g}_2(x)=&\left(\beta(x)-P_{\beta}(x)\right)\cdot \left(DP_{k+1}(x',x_n-P_{\Omega}(x'))-P_{k+1,n}(x',x_n-P_{\Omega}(x'))DP_{\Omega}(x')\right),\\
\bar{g}_2(x)=&\beta(0)\cdot DP_{k+1}(x',x_n-P_{\Omega}(x')),\\
\hat{g}_2(x)=&\left(P_{\beta}(x)-\beta(0)\right)\cdot \left(DP_{k+1}(x',x_n-P_{\Omega}(x'))-P_{k+1,n}(x',x_n-P_{\Omega}(x'))DP_{\Omega}(x')\right)\\
&-\beta(0)\cdot DP_{\Omega}(x')P_{k+1,n}(x',x_n-P_{\Omega}(x')).
\end{aligned}
\end{equation*}

By \cref{e3.20} and noting $P_{k+1}\in \mathcal{HP}_{k+1}$,
\begin{equation}\label{e3.27}
|\tilde{g}_2(x)|\leq C|x|^{k+l-1+\alpha},~~\forall ~x\in (\partial \Omega)_1.
\end{equation}
Since $\beta(0)\cdot DP_{k+1}(x)$ is a $k$-form and $\partial \Omega\in C^{l,\alpha}(0)$,
\begin{equation}\label{e3.22}
|\bar{g}_2(x)|=|\beta(0)\cdot DP_{k+1}(x',x_n-P_{\Omega})|\leq C|x|^{k+l-1+\alpha},~~\forall ~x\in (\partial \Omega)_1.
\end{equation}
Next, since $\hat{g}$ is a polynomial, it can be verified easily that
\begin{equation}\label{e3.26}
\hat{g}_2(0)=\cdots=|D^{k}\hat{g}(0)|=0.
\end{equation}
Then by \crefrange{e3.11}{e3.26}, we have $g_2\in C^{k+l-1,\alpha}(0)$ and
\begin{equation*}
g_2(0)=\cdots=|D^{k}g_2(0)|=0.
\end{equation*}

By \Cref{le-3.1} again, $u_2\in C^{k+2,\alpha}(0)$. That is, there exists $P_{k+2}\in \mathcal{HP}_{k+2}$ such that
\begin{equation*}
  |u_2(x)-P_{k+2}(x)|\leq C|x|^{k+2+\alpha},~~~\forall ~x\in (\Omega)_1
\end{equation*}
and
\begin{equation}\label{e3.23}
\beta(0)\cdot DP_{k+2}(x',0)\equiv \mathbf{\Pi}_{k+1}(g_2(x',0)).
\end{equation}
I.e., $(\beta(0)\cdot DP_{k+2}-\mathbf{\Pi}_{k+1}(g_2))(x)$ is a $(k+1)$-form.

Set
\begin{equation*}
\begin{aligned}
u_3(x)&=u_2-P_{k+2}(x',x_n-P_{\Omega})\\
&=u_1(x)-P_{k+1}(x',x_n-P_{\Omega})-P_{k+2}(x',x_n-P_{\Omega}).
\end{aligned}
\end{equation*}
Hence, $u_3(0)=\cdots=|D^{k+2}u_3(0)|=0$. In addition, $u_3$ is a viscosity solution of
\begin{equation*}
\left\{\begin{aligned}
&\Delta u_3=f_3&& ~~\mbox{in}~~\Omega_1;\\
&\beta\cdot Du_3=g_3&& ~~\mbox{on}~~(\partial \Omega)_1,
\end{aligned}\right.
\end{equation*}
where $f_3\in C^{k+l-2,\alpha}(0)$,
\begin{equation*}
  \begin{aligned}
g_3(x)=&g_2(x)-\beta(x)\cdot\left(DP_{k+2}(x',x_n-P_{\Omega}(x'))
-P_{k+2,n}(x',x_n-P_{\Omega}(x'))DP_{\Omega}(x')\right)\\
=&g_2(x)-\tilde{g}_3(x)-\bar{g}_3(x)-\hat{g}_3(x),\\
\end{aligned}
\end{equation*}
and
\begin{equation*}
\begin{aligned}
\tilde{g}_3(x)=&\left(\beta(x)-P_{\beta}(x)\right)\cdot \left(DP_{k+2}(x',x_n-P_{\Omega}(x'))-P_{k+2,n}(x',x_n-P_{\Omega}(x'))DP_{\Omega}(x')\right),\\
\bar{g}_3(x)=&\beta(0)\cdot DP_{k+2}(x',x_n-P_{\Omega}(x')),\\
\hat{g}_3(x)=&\left(P_{\beta}(x)-\beta(0)\right)\cdot \left(DP_{k+2}(x',x_n-P_{\Omega}(x'))-P_{k+2,n}(x',x_n-P_{\Omega}(x'))DP_{\Omega}(x')\right)\\
&-\beta(0)\cdot DP_{\Omega}(x')P_{k+2,n}(x',x_n-P_{\Omega}(x')).
\end{aligned}
\end{equation*}

Similar to the above argument,
\begin{equation*}
|\tilde{g}_3(x)|\leq C|x|^{k+l-1+\alpha},~~\forall ~x\in (\partial \Omega)_1,
\end{equation*}
and
\begin{equation*}
\hat{g}_3(0)=\cdots=|D^{k+1}\hat{g}(0)|=0.
\end{equation*}
In addition, by \cref{e3.23}, for any $x\in (\partial \Omega)_1$,
\begin{equation*}
\begin{aligned}
  |\bar{g}_3(x)-\mathbf{\Pi}_{k+1}(g_2)(x',x_n-P_{\Omega}(x'))|
&\leq C|x|^{k+l+\alpha}.
\end{aligned}
\end{equation*}
Hence,
\begin{equation*}
  \begin{aligned}
    |g_3(x)|=&|g_2(x)-\tilde{g}_3(x)-\bar{g}_3(x)-\hat{g}_3(x)|\\
=&|g_2(x)-\mathbf{\Pi}_{k+1}(g_2)(x)+\mathbf{\Pi}_{k+1}(g_2)(x',P_{\Omega}(x'))
    +\mathbf{\Pi}_{k+1}(g_2)(x',x_n-P_{\Omega}(x'))\\
    &-\tilde{g}_3(x)-\bar{g}_3(x)-\hat{g}_3(x)|\\
    \leq &C|x|^{k+2}.
  \end{aligned}
\end{equation*}
That is, $g_3(0)=\cdots=|D^{k+1}g_3(0)|=0.$ By virtue of \Cref{l-2.2} again, $u_3\in C^{k+3,\alpha}(0)$. By similar arguments again and again, $u\in C^{k+l,\alpha}(0)$ eventually and \cref{e1.3} holds. Therefore, the proof of \Cref{th1.3} is completed.\qed~\\

The \Cref{th1.5} is an easy consequence of \Cref{th1.3}.~\\
\noindent\textbf{Proof of \Cref{th1.5}.} In fact, we can always assume that $u(0)=|Du(0)|=0$ since the boundary condition is a first order equation. Let $\tilde{u}(x)=u(x)-u(0)-Du(0)\cdot x$ and then $\tilde{u}$ satisfies
\begin{equation*}
 \left\{ \begin{aligned}
&\Delta \tilde{u}=f ~~&&\mbox{in}~~\Omega\cap B_1;\\
&\beta\cdot D\tilde{u}=\tilde{g} ~~&&\mbox{on}~~\partial \Omega\cap B_1,
 \end{aligned}\right.
\end{equation*}
where $\tilde{g}=g-\beta\cdot Du(0)$. Note that $\beta(0)\cdot Du(0)=g(0)$ and hence
\begin{equation*}
  \tilde{g}(x)=g(x)-g(0)-(\beta(x)-\beta(0))\cdot Du(0).
\end{equation*}
Thus, $\tilde{u}(0)=|D\tilde{u}(0)|=\tilde{g}(0)=0$ and $\tilde{g}\in C^{k-1,\alpha}(0)$. By \Cref{th1.3} with $l=1$, we arrive at the conclusion of \Cref{th1.5}.~\qed~\\

Finally, as an application to the regularity of free boundaries in one phase problems, we prove \Cref{th1.4}.~\\
\noindent\textbf{Proof of \Cref{th1.4}.} Since $u=0$ on $\partial \Omega\cap B_1$ and $\partial \Omega\in C^{1,\alpha}$, by the boundary pointwise regularity for Dirichlet problems, $u\in C^{1,\alpha}(\bar{\Omega}\cap B_1)$. Hence, we have $|Du|=1$ on $(\partial \Omega)_1$ in the classical sense. Let $v(x)=u(x)-x_n$ and then
\begin{equation}\label{e3.24}
v_n=u_n-1=\left(1-\sum_{i=1}^{n-1}u_i^2\right)^{1/2}-1
=\frac{-\sum_{i=1}^{n-1}u_i^2}{1+\left(1-\sum_{i=1}^{n-1}u_i^2\right)^{1/2}}~~\mbox{ on }~(\partial \Omega)_1.
\end{equation}

Since $u\in C^{1,\alpha}(0)$, $u_i\in C^{\alpha}(0)$ and hence $u_i^2\in C^{2\alpha}(0)$. By the boundary pointwise regularity for oblique derivative problems (see \Cref{th1.3-1}), $v\in C^{1, 2\alpha}(0)$ and then $u\in C^{1, 2\alpha}(0)$. Similarly, $u\in C^{1, 2\alpha}(x_0)$ for any $x_0\in (\partial \Omega)_1$. Note that $u_n(0)=1$ and then $u_i/u_n\in C^{2\alpha}(\bar\Omega \cap B_r)$ for some $r>0$. Thus, $\partial \Omega\cap B_r\in C^{1,2\alpha}$. Likewise, $(\partial \Omega)_1\in C^{1,2\alpha}$.

Repeat above argument. Note $u_i^2\in C^{4\alpha}(0)$. By \Cref{th1.3-1}, $v\in C^{1, 4\alpha}(0)$ and we have $(\partial \Omega)_1\in C^{1,4\alpha}$. After finite steps, $u\in C^{2,\tilde\alpha}$ and $(\partial \Omega)_1\in C^{2,\tilde\alpha}$ for some $0<\tilde\alpha<1$.

Note that $u_i(0)=0$ for $1\leq i\leq n-1$. Then $u_i^2\in C^{2,\tilde\alpha}(0)$. By \Cref{th1.5}, $v\in C^{3,\tilde\alpha}(0)$. Hence $u\in C^{3,\tilde\alpha}(\Omega_1)$ and $\partial \Omega\in C^{3,\tilde\alpha}$. Thus, $u_i^2\in C^{3,\tilde\alpha}$ then $u\in C^{4,\tilde\alpha}$ and $\partial \Omega\in C^{4,\tilde\alpha}$. By iteration arguments, we have $u\in C^\infty$ and $\partial \Omega\in C^{\infty}$ eventually.~\qed~\\

\begin{remark}\label{re3.7}
Since $u=0$ and $|Du|=1$ on $(\partial \Omega)_1$, $\partial u/\partial \nu=1$ on $(\partial \Omega)_1$ where $\nu$ is the inner normal. Thus, maybe a more natural idea of proving \Cref{th1.4} is to consider the Neumann problem:
\begin{equation*}
\left\{\begin{aligned}
&   \Delta u=0~~~~\mbox{in}~~\Omega \cap B_1;\\
&   \frac{\partial u}{\partial \nu}=1~~~~\mbox{on}~~\partial \Omega\cap B_1.
\end{aligned}\right.
\end{equation*}
However, the smoothness of $\nu$ depends on the smoothness of $(\partial \Omega)_1$. Hence, we can't improve the regularity of solutions and $(\partial \Omega)_1$. Instead, $\beta\equiv e_n$ if considering \cref{e3.24}.
\end{remark}
~\\

\noindent\textbf{Data availability statement} Data sharing not applicable to this article as no datasets were generated or analysed during the current study.

\bibliographystyle{amsplain}
\bibliography{PDE}

\providecommand{\bysame}{\leavevmode\hbox to3em{\hrulefill}\thinspace}
\providecommand{\MR}{\relax\ifhmode\unskip\space\fi MR }
\providecommand{\MRhref}[2]{%
  \href{http://www.ams.org/mathscinet-getitem?mr=#1}{#2}
}
\providecommand{\href}[2]{#2}
\begin{thebibliography}{10}

\bibitem{MR0125307}
S.~Agmon, A.~Douglis, and L.~Nirenberg, \emph{Estimates near the boundary for
  solutions of elliptic partial differential equations satisfying general
  boundary conditions. {I}}, Comm. Pure Appl. Math. \textbf{12} (1959),
  623--727. \MR{0125307}

\bibitem{MR803243}
Ioannis Athanasopoulos and Luis~A. Caffarelli, \emph{A theorem of real analysis
  and its application to free boundary problems}, Comm. Pure Appl. Math.
  \textbf{38} (1985), no.~5, 499--502. \MR{803243}

\bibitem{MR1005611}
Luis~A. Caffarelli, \emph{Interior a priori estimates for solutions of fully
  nonlinear equations}, Ann. of Math. (2) \textbf{130} (1989), no.~1, 189--213.
  \MR{1005611}

\bibitem{MR1351007}
Luis~A. Caffarelli and Xavier Cabr\'e, \emph{Fully nonlinear elliptic
  equations}, American Mathematical Society Colloquium Publications, vol.~43,
  American Mathematical Society, Providence, RI, 1995. \MR{1351007}

\bibitem{MR1118699}
Michael~G. Crandall, Hitoshi Ishii, and Pierre-Louis Lions, \emph{User's guide
  to viscosity solutions of second order partial differential equations}, Bull.
  Amer. Math. Soc. (N.S.) \textbf{27} (1992), no.~1, 1--67. \MR{1118699}

\bibitem{MR4093736}
D.~De~Silva and O.~Savin, \emph{A short proof of boundary {H}arnack principle},
  J. Differential Equations \textbf{269} (2020), no.~3, 2419--2429.
  \MR{4093736}

\bibitem{MR3393271}
Daniela De~Silva and Ovidiu Savin, \emph{A note on higher regularity boundary
  {H}arnack inequality}, Discrete Contin. Dyn. Syst. \textbf{35} (2015),
  no.~12, 6155--6163. \MR{3393271}

\bibitem{MR3198649}
Alessio Figalli and Henrik Shahgholian, \emph{A general class of free boundary
  problems for fully nonlinear elliptic equations}, Arch. Ration. Mech. Anal.
  \textbf{213} (2014), no.~1, 269--286. \MR{3198649}

\bibitem{MR0440187}
D.~Kinderlehrer and L.~Nirenberg, \emph{Regularity in free boundary problems},
  Ann. Scuola Norm. Sup. Pisa Cl. Sci. (4) \textbf{4} (1977), no.~2, 373--391.
  \MR{0440187}

\bibitem{MR3780142}
Dongsheng Li and Kai Zhang, \emph{Regularity for fully nonlinear elliptic
  equations with oblique boundary conditions}, Arch. Ration. Mech. Anal.
  \textbf{228} (2018), no.~3, 923--967. \MR{3780142}

\bibitem{lian2020pointwise}
Yuanyuan Lian, Lihe Wang, and Kai Zhang, \emph{Pointwise regularity for fully
  nonlinear elliptic equations in general forms, arxiv: 2012.00324}, 2020.

\bibitem{MR4088470}
Yuanyuan Lian and Kai Zhang, \emph{Boundary pointwise {$C ^{1,\alpha}$} and
  {$C^{2,\alpha}$} regularity for fully nonlinear elliptic equations}, J.
  Differential Equations \textbf{269} (2020), no.~2, 1172--1191. \MR{4088470}

\bibitem{MR906819}
Gary~M. Lieberman, \emph{Local estimates for subsolutions and supersolutions of
  oblique derivative problems for general second order elliptic equations},
  Trans. Amer. Math. Soc. \textbf{304} (1987), no.~1, 343--353. \MR{906819}

\bibitem{MR2254613}
Emmanouil Milakis and Luis~E Silvestre, \emph{{Regularity for fully nonlinear
  elliptic equations with {N}eumann boundary data}}, Comm. Partial Differential
  Equations \textbf{31} (2006), no.~7-9, 1227--1252.

\bibitem{MR3980853}
Gabrielle Nornberg, \emph{{$C^{1,\alpha}$} regularity for fully nonlinear
  elliptic equations with superlinear growth in the gradient}, J. Math. Pures
  Appl. (9) \textbf{128} (2019), 297--329. \MR{3980853}

\bibitem{MR2962060}
Arshak Petrosyan, Henrik Shahgholian, and Nina Uraltseva, \emph{Regularity of
  free boundaries in obstacle-type problems}, Graduate Studies in Mathematics,
  vol. 136, American Mathematical Society, Providence, RI, 2012. \MR{2962060}

\bibitem{MR0478079}
Aleksey~Vasil'yevich Pogorelov, \emph{The {M}inkowski multidimensional
  problem}, V. H. Winston \& Sons, Washington, D.C.; Halsted Press [John Wiley
  \&\ Sons], New York-Toronto-London, 1978, Translated from the Russian by
  Vladimir Oliker, Introduction by Louis Nirenberg, Scripta Series in
  Mathematics. \MR{0478079}

\bibitem{MR2983006}
O.~Savin, \emph{Pointwise {$C^{2,\alpha}$} estimates at the boundary for the
  {M}onge-{A}mp\`ere equation}, J. Amer. Math. Soc. \textbf{26} (2013), no.~1,
  63--99. \MR{2983006}

\bibitem{MR2334822}
Ovidiu Savin, \emph{Small perturbation solutions for elliptic equations}, Comm.
  Partial Differential Equations \textbf{32} (2007), no.~4-6, 557--578.
  \MR{2334822}

\bibitem{MR3246039}
Luis Silvestre and Boyan Sirakov, \emph{Boundary regularity for viscosity
  solutions of fully nonlinear elliptic equations}, Comm. Partial Differential
  Equations \textbf{39} (2014), no.~9, 1694--1717. \MR{3246039}

\bibitem{MR2592289}
Boyan Sirakov, \emph{Solvability of uniformly elliptic fully nonlinear {PDE}},
  Arch. Ration. Mech. Anal. \textbf{195} (2010), no.~2, 579--607. \MR{2592289}

\bibitem{MR1135923}
Lihe Wang, \emph{On the regularity theory of fully nonlinear parabolic
  equations. {I}}, Comm. Pure Appl. Math. \textbf{45} (1992), no.~1, 27--76.
  \MR{1135923}

\bibitem{MR1139064}
\bysame, \emph{On the regularity theory of fully nonlinear parabolic equations.
  {II}}, Comm. Pure Appl. Math. \textbf{45} (1992), no.~2, 141--178.
  \MR{1139064}

\bibitem{MR1151267}
\bysame, \emph{On the regularity theory of fully nonlinear parabolic equations.
  {III}}, Comm. Pure Appl. Math. \textbf{45} (1992), no.~3, 255--262.
  \MR{1151267}

\end{thebibliography}

\end{document}